\documentclass[11pt]{article}
\usepackage{amsfonts,latexsym}
\usepackage{theorem}

\newcommand{\Section}[1]{
\section{#1}}
\newtheorem{thm}{Theorem}
\newtheorem{prop}[thm]{Proposition}

\newtheorem{lemma}[thm]{Lemma}
{\theorembodyfont{\normalfont\rmfamily}
\newtheorem{definition}[thm]{Definition}
\newtheorem{remark}[thm]{Remark}

}
\makeatletter
\def\ack{\vspace{.5\baselineskip}\noindent{\theorem@headerfont
Acknowledgement}\ \ }
\newenvironment{proof}[1][]%
{\def\proof@temp{#1}\par\noindent
\textsc{Proof}\ifx\proof@temp\@empty\else\ ({\proof@temp})\fi\hspace{1em}}
{\hphantom{xxx}\hfill~ {$\Box$}\par\vspace{.4\baselineskip}}
\newcounter{Item}

\newcommand{\email}[1]{{\small\tt #1}}

\def\operatorname#1{\mathop{\operator@font #1}\nolimits}%
\makeatother

\newcommand{\A}{\mathcal{A}}
\renewcommand{\H}{\mathcal{H}}
\newcommand{\E}{\mathcal{E}}
\newcommand{\D}{\mathcal{D}}
\newcommand{\W}{\mathcal{W}}
\newcommand{\CALS}{\mathcal{S}}
\newcommand{\B}{\mathcal{B}}
\renewcommand{\L}{\mathcal{L}}
\newcommand{\C}{\mathbb{C}}
\newcommand{\R}{\mathbb{R}}
\newcommand{\Z}{\mathbb{Z}}
\newcommand{\PO}{P}
\newcommand{\ehat}[1]{\hat{E}({L^{#1}})}
\newcommand{\CL}{{\mathcal C}_L}
\newcommand{\Ci}[1]{C^{\infty}(#1)}

\newcommand{\ad}{\operatorname{ad}}

\newcommand{\Id}{\operatorname{Id}}

\newcommand{\cech}{\v Cech}
\newcommand{\CiM}{C^\infty(M)}
\newcommand{\bbnu}{[\![\nu]\!]}
\newcommand{\coh}[1]{H^{#1}(M;\R)}

\newcommand{\alphabeta}{{\alpha\beta}}
\newcommand{\alphagamma}{{\alpha\gamma}}
\newcommand{\betaalpha}{{\beta\alpha}}
\newcommand{\gammabeta}{{\gamma\beta}}
\newcommand{\alphabetagamma}{{\alpha\beta\gamma}}
\newcommand{\gammabetaalpha}{{\gamma\beta\alpha}}

\newcommand{\cob}{\partial}
\newcommand{\mdot}{\mathop{\raise.4ex\hbox{.}}}

\newcommand{\AkB}{\widehat{A}*_k \widehat{B}}
\newcommand{\AB}{\widehat{A} * \widehat{B}}
\newcommand{\G}{G_0}
\newcommand{\K}{K_0}
\newcommand{\GC}{G}
\newcommand{\KC}{K}
\newcommand{\g}{\mathfrak{g}}
\renewcommand{\k}{\mathfrak{k}}

\newcommand{\gc}{\mathfrak{g}}
\newcommand{\kc}{\mathfrak{k}}
\newcommand{\m}{{\mathfrak{m}}}
\newcommand{\epsilonk}{\epsilon^{(k)}}
\newcommand{\hatf}{\widehat f\,}
\newcommand{\half}{\raise0.35ex\hbox{$\scriptstyle{1\over2}$}}
\newcommand{\quarter}{\raise0.35ex\hbox{$\scriptstyle{1\over4}$}}

\newdimen{\addresswidth}
\title{ Variations on \\
 deformation quantization }

\author{
S. Gutt
\thanks{Research supported by the Communaut\'e fran\c caise de
Belgique, through an Action de Recherche Concert\'ee de la Direction de
la Recherche Scientifique.}\\
\small \hbox{\parbox[t]{1.9in}{\begin{center}%
Universit\'e Libre de Bruxelles\\[3pt]
Campus Plaine, CP 218\\[3pt]
bd du Triomphe\\[3pt]
1050 Brussels, Belgium\\[3pt]
\email{sgutt\char64ulb.ac.be}\end{center}}}
\hbox{\parbox[t]{.6in}{\begin{center}\rm and\end{center}}}
\hbox{\parbox[t]{1.9in}{\begin{center}%
Universit\'e de Metz\\[3pt]
Ile du Saulcy\\[3pt]
57045 Metz Cedex 01, France\end{center}}}
}

\date{}

\begin{document}
{\renewcommand{\baselinestretch}{1} 
\maketitle\thispagestyle{empty}

\begin{abstract}
I was asked  by the organisers to present some 
 aspects of Deformation Quantization. Mosh\'e 
 has pursued, for more than 25 years, a research program based
 on the idea that physics progresses in stages, and one goes from one
 level of the theory to the next one by a deformation,
  in the mathematical sense of the word, to be defined in an
  appropriate category. His study of deformation theory applied
 to mechanics started in 1974 and led to spectacular developments
 with the deformation quantization programme.
  
 \medskip
  I first met Mosh\'e at a  conference
  in  Li\`ege in 1977. A few months later he became my thesis ``codirecteur''.
  Since then he has been one of my closest friends, present at all stages
  of my personal and mathematical life. I miss him....

  \medskip  
I have chosen, in this presentation of Deformation Quantization,
to focus on 3 points: the uniqueness --up to equivalence--
of a universal  star product (universal in the sense of Kontsevich) 
on the dual of a Lie algebra, the cohomology classes introduced by
Deligne for equivalence classes of differential star products on a
symplectic manifold and the construction of some convergent
star products on Hermitian symmetric spaces. 
Those subjects will appear in a  promenade through the history of 
 existence and equivalence in deformation quantization.
\end{abstract}

\bigskip

\centerline{Moshe Flato Conference, Dijon, September 1999}

\newpage


\Section{Introduction}\label{sect:intro}

Quantization of a classical system is a way to pass from classical
to quantum results.

Classical mechanics, in its  Hamiltonian formulation on
the motion space, has for framework
a symplectic manifold (or more generally a Poisson manifold).
Observables are families of smooth functions on that manifold $M$.
The dynamics is defined in terms of  a Hamiltonian $H \in C^\infty(M)$
and the time evolution
of an observable $f_t \in C^\infty(M\times\R)$
is governed by the equation :
$$
{{d}\over{dt}}f_t = -\left\{H,f_t\right\}.
$$

Quantum mechanics, in its usual Heisenberg's formulation, has for 
framework a Hilbert space (states are rays in that space).
 Observables are 
families   of selfadjoint operators on the Hilbert space.
The dynamics is defined in terms of a Hamiltonian $H$, which 
is a selfadjoint operator, and the time evolution of an
observable $A_t$ is governed by the equation~:
$$
{{dA_t}\over{dt}} = {{i}\over{\hbar}} [H,A_t].
$$

A natural suggestion for quantization is a correspondence
$Q\colon f\mapsto Q(f)$ mapping a function $f$ to a self adjoint
operator $Q(f)$ on a Hilbert space $\H$ in such a way that
$Q(1)=\Id$ and
\[
[Q(f),Q(g)]=i\hbar Q(\{f,g\}).
\]
There is no such correspondence defined
on all smooth functions on $M$ when one puts  an irreducibility
requirement which is necessary not to violate Heisenberg's principle.

Different mathematical treatments of quantization appeared to 
deal with this problem:
\begin{itemize}
\item Geometric Quantization of Kostant and Souriau. This proceeds
in two steps; first prequantization of a symplectic manifold
$(M,\omega)$ where one builds a Hilbert space and a correspondence
$Q$ as above defined on all smooth functions on $M$ but with
no irreducibility, then polarization to ``cut down the number of variables''.
One succeeds to quantize only a small class of functions.

\item Berezin's quantization where one builds on a particular class of
K\"ahler manifolds a family of associative algebras
using a symbolic calculus, i.e. a dequantization procedure.

 \item Deformation Quantization introduced by Flato,  Lichnerowicz
 and  Sternheimer in \cite{FLS} and developed in \cite{BayenB} where they 

`` suggest that quantization be understood as a deformation of the
structure of the algebra of classical observables rather than a radical change
in the nature of the observables.''

\end {itemize}
This deformation approach to quantization is 
part of a  general deformation approach to physics. This was
one of the seminal ideas stressed by Moshe: one looks
at some level of a theory in physics as a deformation of another
level \cite{Flato}.

Deformation quantization is  defined in terms
of a star product which is a formal
deformation of the algebraic structure of the space of smooth functions
on a Poisson manifold. The associative structure given by the
usual product of functions and the Lie structure given by the Poisson
bracket are simultaneously deformed. 

The plan of this presentation is the following :  

\begin{itemize}
\item  Examples and existence of star products. After some definitions,
I recall the history of the proofs of existence and through the
history of some  constructions
 of star products, starting from the Moyal star product on a vector
 space endowed with a constant Poisson structure. 
  I give the standard star product on the dual
of a Lie algebra, and recall Kontsevich star product for any Poisson structure
on a vector space. I prove by elementary methods
that all universal star products on the dual of a Lie algebra are
essentially equivalent.

\item Equivalence of star products. I first wander through the history of
the parametrization of equivalence classes of star products. I define
Deligne's cohomology classes associated to differential star products on
symplectic manifolds and show by \cech ~ methods how this yields an intrinsic
parametrization of the equivalence classes of differential star products
 and how this allows to study automorphisms of a star product
 (in the symplectic framework). Finally, I define
  a generalised moment map for star products.

\item Convergence of star products. I study mostly  the convergence of 
a Berezin type star product on Hermitian symmetric spaces. The construction
of such a star product involves a correspondence between operators (on
the Hilbert space given by geometric quantization) and functions (their
Berezin symbols).  A parameter is introduced in the construction
(generalising the power of the line bundle in geometric quantization). 
Asymptotic expansion in this parameter on a large algebra of
functions yields a deformed product. One proves associativity  and
 convergence of this product.
\end{itemize}
I would like to thank my friends Georges Pinczon
and Daniel Sternheimer who brought many improvements to this presentation.

\section{  Examples and existence of star products}\label{sect:ex}

\begin{definition}
A \textbf{Poisson bracket} defined on the space of
smooth functions on a manifold $M$, 
is a $\R$- bilinear map on $\CiM$,
$(u,v)\mapsto \{u,v\}$ such that for any $u,v,w \in \CiM$:

- $\{u,v\}=-\{v,u\}$;

- $\{\{u,v\},w\}+\{\{v,w\},u\}+\{\{w,u\},v\}=0$;

- $\{u,vw\}=\{u,v\}w+\{u,w\}v$.

\noindent A Poisson bracket is given in terms of a contravariant skew 
symmetric 2-tensor $\PO$ on $M$, called the \textbf{Poisson tensor}, by
$$
\{u,v\} = \PO(du\wedge dv).
$$
The Jacobi identity for the Poisson bracket Lie algebra is equivalent
to the vanishing of  the Schouten bracket :
$$
[\PO,\PO]=0.
$$
(The Schouten bracket is the extension -as a graded derivation for the 
exterior product- of the bracket of vector fields to skewsymmetric 
contravariant tensor fields.) 

\noindent A \textbf{Poisson manifold}, denoted $(M,\PO)$, is a manifold $M$
with a Poisson bracket defined by the Poisson tensor $\PO$. 
\end{definition}

\medskip
A particular class of Poisson manifolds, essential in classical mechanics,
 is the class
 of {\bf symplectic manifolds}. If $(M,\omega)$ is a symplectic manifold
 (i.e.  $\omega$ is a closed nondegenerate $2$-form on $M$) and 
 if $u,v \in \CiM$,  the Poisson bracket of $u$ and $v$ is defined by
$$
\{u,v\} : = X_u(v) = \omega(X_v,X_u),
$$
where $X_u$ denotes the  Hamiltonian vector field 
corresponding to the function $u$,
i.e. such that $
i(X_u)\omega = du.$
In coordinates the components  of the corresponding Poisson tensor 
$\PO^{ij}$ form the inverse matrix of
the components $\omega_{ij}$ of $\omega$.

\medskip
{\bf Duals of Lie algebras} form  the class of linear Poisson manifolds.
If $\g$ is a Lie algebra then its dual ${\g}^*$ is endowed with
the Poisson tensor $\PO$ defined by 
$$
\PO_\xi(X,Y) : =\xi([X,Y])
$$
for $X,Y\in\g=(\g^*)^*\sim(T_\xi\g^*)^*$.

\begin{definition}\label{def:star} (Bayen et al.\ \cite{BayenB})\ \
A \textbf{star product} on $(M,\PO)$ is a bilinear map
\[
N\times N \to N\bbnu, \qquad (u,v) \mapsto u*v = u*_\nu v : =
\sum_{r\ge0} \nu^rC_r(u,v)\]
where $N=\CiM$, such that
\begin{enumerate}
\item when the map is extended $\nu$-linearly (and 
 continuously in the $\nu$-adic topology) to $N\bbnu\times N\bbnu$
 it is formally associative:
\[
(u*v)*w = u*(v*w);
\]
\item (a)\ \ $C_0(u,v) = uv$,\ \ \ (b)\ \ $C_1(u,v)-C_1(v,u) = \{u,v\}$;
\item $1*u = u*1 = u$.
\end{enumerate}
When the $C_r$'s are bidifferential operators on $M$, one speaks  of
 a {\bf differential star product}.  

\end{definition}

\begin{remark}
 A star product
can also be defined not on the whole of $\CiM$ but on any subspace $N$ of
it which is stable under pointwize multiplication and Poisson bracket.

Requiring  differentiability of the cochains is essentially the same as
requiring them to be local \cite{CDG}.

In (b) we follow Deligne's normalisation for $C_1$: its
skew symmetric part is $\frac12 \{\,,\,\}$. In the original definition it
was equal to the Poisson bracket. One finds in the literature other
normalisations such as $\frac{i}{2} \{\,,\,\}$. All these amount to
a rescaling of the parameter $\nu$.

One assumed also the parity condition
$C_r(u,v)=(-1)^r C_r(v,u)$ in the earliest definition.

Property (b) above implies that the centre of $\CiM\bbnu$, when the latter
is viewed as an algebra with multiplication $*$, is a series whose terms
Poisson commute with all functions, so is an element of $\R\bbnu$ when
$M$ is symplectic and connected.

Properties (a) and (b) of Definition \ref{def:star} imply
 that the \textbf{star
commutator} defined by $[u,v]_* = u*v - v*u $,
 which obviously makes $\CiM\bbnu$ into a Lie algebra, has the form $
[u,v]_* = \nu \{u,v\} + \dots $
so that repeated bracketing leads to higher and higher order terms. This
makes $\CiM\bbnu$ an example of a \textbf{pronilpotent Lie algebra}. 
We denote the \textbf{star adjoint
representation} $ad_* u \,\, (v) = [u,v]_*.$
\end{remark}

\subsection{The Moyal star product on $\R^n$}

The simplest example of a deformation quantization
is the Moyal product for the Poisson structure $\PO$ on a vector space
$V=\R^m$ with constant coefficients:
 \[
 \PO=\sum_{i,j} \PO^{ij} {\partial}_{i} \wedge
 {\partial}_{j},\,\,\,\PO^{ij}=-\PO^{ji}\in \R
\] 
 where ${\partial}_i={\partial}/{\partial} x^i$ is the partial derivative 
 in the direction of the coordinate $x^i,\,\,\,i=1,\dots, n$.
 The formula for {\bf the Moyal product} is
 \begin{equation}
 (u*_Mv)(z)=
 \left.\exp \left(\frac{\nu}{2} \PO^{rs}
 \partial_{x^r}\partial_{y^s}\right)(u(x)v(y))\right\vert_{x=y=z}. 
\end{equation}

When $\PO$ is non degenerate (so $V=\R^{2n}$),  the space of formal power
 series of polynomials on $V$ with Moyal product
is called the Weyl algebra $W=(S(V)[[\nu]],*_M)$.

This example  comes from the composition of operators via Weyl's quantization.
Weyl's correspondence associates to certain functions $f$ on $\R^{2n}$
an operator $W(f)$ on $L^2(\R^n)$
\[
W(f)=\int\tilde{f}(\xi,\eta) \exp\left(\frac{i(\xi P+\eta Q)}{\hbar}\right)d\xi~d\eta
\]
where $\tilde{f}$ is the inverse Fourier transform of $f$.
Then
\[
W(f)\circ W(g)=W(f*_Mg)  \quad (\nu={i\hbar}).
\]
In fact, Moyal had  used in 1949 the deformed bracket which corresponds to the commutator
of operators to  study quantum statistical mechanics. The Moyal product first appeared 
in Groenewold.
\bigskip

$\bullet$ In 1974, Flato, Lichnerowicz and Sternheimer \cite{FLS1}  
studied deformations
of the Lie algebra structure defined by the Poisson bracket on the algebra
$N$
of smooth functions on a symplectic manifold; they studied
1-differential deformations because  the relevant cohomology,
i.e. the 1-differential Chevalley cohomology of the Lie algebra $N$
 with values in $N$ for the adjoint representation, was known \cite{Lich}.  

In 1975, Vey \cite{Vey} pursued their work in the differential context:
 he  constructed a differential deformation on $M=\R^{2n}$
 which turns out to be the
 Moyal bracket $\{u,v\}_M:=\frac1{\nu}(u*_Mv-v*_Mu)$ and he proved
 that there exists a differential deformation on
 a symplectic manifold when its third de Rham cohomology
space is trivial ($b_3(M):=\dim H^3(M;\R)=0$). He also 
reconstructed the Moyal product on $\R^{2n}$.

This  opened the path to  deformation quantization presented 
in 1976 by Flato, Lichnerowicz and Sternheimer \cite{FLS}.

\medskip
$\bullet$  In 1978, in their seminal paper about deformation 
quantization  \cite{BayenB},
 Bayen, Flato, Fronsdal, Lichnerowicz and Sternheimer proved that
Moyal star product can be defined on any symplectic manifold $(M, \omega)$
which admits a  symplectic connection $\nabla$ (i.e. a linear connection
such that 
$\nabla\omega=0$ and
the torsion of $\nabla$ vanishes) with no curvature. They also built star
products on some quotient of $\R^{2n}$.

\medskip
$\bullet$  In 1979, Neroslavsky and Vlassov \cite{NeroVla} proved 
with Lichnerowicz that on any symplectic manifold
with $b_3(M)=0$, there exists a differential star product.

\subsection{Hochschild cohomology}

The study of star products  on a manifold $M$ used Gerstenhaber theory of 
deformations  \cite{Gerst} of associative algebras. 
This uses the Hochschild cohomology of the algebra, here
$\CiM$ with values in $\CiM$, where $p$-cochains are $p$-linear maps from
$(\CiM)^p$ to $\CiM$ and where the {\bf Hochschild
coboundary operator} maps the $p$-cochain $C$ to the $p+1$-cochain
\begin{eqnarray*}
&&(\cob C)(u_0,\dots,u_p) = u_0 C(u_1,\dots,u_p) \\
 &&~+ \sum_{r=1}^p
(-1)^r C(u_0,\dots,u_{r-1}u_r,\dots,u_p)
+ (-1)^{p + 1} C(u_0,\dots,u_{p-1})u_p.
\end{eqnarray*}
For differential star products, we consider differential cochains, i.e. 
given by differential operators on each argument.
The associativity condition for a star product at order $k$ in the
parameter $\nu$ reads
\[
(\cob C_k)(u,v,w)=\sum_{r+s=k,r,s>0}\left(~
C_r(C_s(u,v),w)-C_r(u,C_s(v,w))~\right).
\]
If one has cochains $C_j,j<k$ such that the star product they define
is associative to order $k-1$, then the right hand
side above is a cocycle ($\cob$(RHS)$=0$) and one can extend the star
product to order $k$ if it is a coboundary (RHS$=\cob(C_k))$.

Denoting by $m$ the usual multiplication of functions,
and writing $*=m+C$ where $C$ is a  formal series of 
multidifferential operators
($C\in D_{poly}(M)[[\nu]]$)  the associativity also reads  
$\cob C=[C,C]$ where the
bracket on the right hand side is the graded Lie algebra bracket on
$D_{poly}(M)[[\nu]]$.

\begin{thm} (Vey \cite{Vey})\label{cohomoVey}
Every differential $p$-cocycle $C$ on a manifold $M$ is the sum of the
coboundary of a differential (p-1)-cochain and a $1$-differential
skewsymmetric $p$-cocycle $A$:
$$
C = \cob B + A.
$$
In particular, a cocycle is a coboundary if and only if its total
skewsymmetrization, which is automatically $1$-differential in each
argument, vanishes. Also 
$$
H^p_\mathrm{diff}(\CiM,\CiM) = \Gamma(\Lambda^p TM).
$$
Furthermore (\cite{CaGulocsym}),given a connection $\nabla$ on $M$,
$B$ can be defined from $C$ by universal formulas.
\end{thm}

 By universal, we mean the following: any $p$-differential operator $D$
of order maximum $k$ in each argument can be written
\[
D(u_1,\ldots,u_p)=\sum_{\vert\alpha_1\vert<k\ldots\vert\alpha_p\vert<k}
  D_{\vert\alpha_1\vert,\ldots,\vert\alpha_p\vert}^{\alpha_1\ldots\alpha_p}
  \nabla_{\alpha_1}u_1\ldots\nabla_{\alpha_p}u_p
\]
where $\alpha$'s are multiindices, 
$D_{\vert\alpha_1\vert,\ldots,\vert\alpha_p\vert}$
are tensors (symmetric in each of the $p$ groups of indices) and
$\nabla_\alpha u=(\nabla\ldots(\nabla u))(\frac{\partial}{\partial x^{i_1}},
\ldots,\frac{\partial}{\partial x^{i_q}})$ when $\alpha=(i_1,\ldots,i_q)$.
We claim that there is a $B$ such that the tensors defining 
$B$ are universally
defined as linear combinations of the tensors defining $C$, 
universally meaning in
a way which is independent of the form of $C$.
 Note that requiring 
differentiability of the cochains is essentially the same as requiring them to be 
local \cite{CDG}.

\smallskip
(An elementary proof of the above theorem can be found in \cite{GR}.)

\begin{remark}
 Behind  theorem \ref{cohomoVey} above, are the following stronger
 results about Hochschild cohomology:
  \begin{thm} 
Let ${\mathcal A}=\CiM$, let $\mathcal{C(A)}$ be the space of continuous
 cochains and $\mathcal{C}_{diff}(\mathcal{A})$ be the space of differential
 cochains. Then
 \begin{itemize}
\item[1)] $\Gamma(\Lambda^p TM)\subset H^p(\CiM,\CiM)$;
\item[2)] the inclusions $\Gamma(\Lambda^p TM)\subset\mathcal{C}_{diff}(\mathcal{A})
\subset\mathcal{C(A)}$ induce isomorphisms in cohomology.
\end{itemize}
\end{thm}
Point 1 follows from the fact that any cochain which is 1-differential in
each argument is a cocycle and that the skewsymmetric part of a coboundary
always vanishes.
The fact that the inclusion $\Gamma(\Lambda TM)\subset
\mathcal{C}_{diff}(\mathcal{A})$
 induces an isomorphism in cohomology is proven by Vey \cite{Vey}; it gives
 theorem \ref{cohomoVey}.
The general result about continuous cochains is due to Connes \cite{Connes}. 
Another proof of Connes result was given by 
Nadaud in \cite{Nadaud}. In the somewhat pathological case of completely 
general cochains the full cohomology does not seem to be known.
\end{remark}

\medskip
$\bullet$ Some examples of star products were built  using inductively the 
     explicit formulas for coboundaries
      on locally symmetric symplectic manifolds
      in \cite{CaGulocsym}; the cochains are given by universal
      expressions in  the symplectic $2$-form,
     its inverse (i.e.  the Poisson
      tensor) and the curvature tensor of the symmetric symplectic connexion. 

\medskip
$\bullet$ We also showed that assuming a homogeneity condition on the cochains
  proves the existence of a star product on the cotangent bundle
   to any parallelisable 
  manifold \cite{CGparalman} and gave explicit formulas 
  for the cotangent bundle to a Lie group.
  The vertical part of this  gives a deformation quantization of any linear
  Poisson manifold, i.e. any dual of a Lie algebra \cite{Gutg}.

\subsection{The standard *-product on $\g^*$}
Let $\g^*$ be the dual of a Lie algebra $\g$. The algebra of polynomials
on $\g^*$ is identified with the symmetric algebra $S(\g)$. One defines
a new associative law on this algebra by a transfer of the product $\circ$ in
the universal enveloping algebra $U(\g)$, via the bijection between
$S(\g)$ and $U(\g)$ given by the total symmetrization $\sigma$ :
\[
\sigma: S(\g)\rightarrow U(\g) X_1  \dots X_k\mapsto
 \frac{1}{k!}\sum_{\rho \in S_k}
X_{\rho(1)}\circ\dots\circ X_{\rho (k)}.
\]
Then $U(\g)=\oplus_{n\ge 0} U_n$ where $U_n:=\sigma (S^n(\g))$ and we decompose
an element $u\in U(\g)$ accordingly $u=\sum u_n$.
We define  for $P\in S^p(\g)$ and $Q\in S^q(\g)$
\begin{equation}\label{starCBH}
P*Q=\sum_{n\ge 0} (\nu)^n \sigma^{-1}((\sigma(P)\circ\sigma(Q))_{p+q-n}).
\end{equation}
This yields a differential star product on $\g^*$ \cite{Gutg}.
Using Vergne's result on the multiplication in $U(\g)$, this
star  product  is characterised  by
\begin{eqnarray*}
&&X*X_1 \dots X_k= XX_1\dots X_k\\
&&\quad+ \sum_{j=1}^k \frac{(-1)^j}{j!}\nu^j
Bj[[[X,X_{r_1}],\dots],X_{r_j}]X_1\dots\widehat{X_{r_1}}
\dots\widehat{X_{r_j}}
\dots X_k
\end{eqnarray*}
where $B_j$ are the Bernouilli numbers.
This star product  can be written with an integral 
formula (for $\nu=2\pi i$)
\cite{Drinfeld}:
\[
u*v(\xi)=\int_{\g\times\g} \hat{u}(X)\hat{v}(Y)e^{2i\pi  \langle\xi, 
CBH(X,Y)\rangle} dXdY
\]
where $\hat{u}(X)=\int_{\g^*}u(\eta) e^{-2i\pi\langle\eta, X\rangle}$
 and where $CBH$ denotes Campbell-Baker-Hausdorff formula for the product
of elements in the group in a logarithmic chart 
($\exp X\exp Y=\exp CBH(X,Y)\quad \forall X,Y\in\g$).

We  call this {\bf the standard (or CBH) star product} on the dual 
of a Lie algebra.

\medskip
$\bullet$
An important property of this star product is its covariance, 
i.e. the fact that
\[
X*Y-Y*X=\nu [X,Y] \quad \forall X,Y \in \g.
\]
In general, when there is a  Hamiltonian action of a connected 
Lie group $G$ on a manifold $M$,
denoting by $\lambda_X$ the function on $M$ corresponding to $X\in\g$,
 a star
product on $M$ is \textbf{covariant} under $G$ if 
$\lambda_X*\lambda_Y -\lambda_Y*\lambda_X=\nu
\lambda_{\{X,Y\}} \forall X,Y\in\g$.
Arnal, Cortet, Molin and Pinczon   have proven in \cite{ACMP}
that for any covariant star product there is a representation of $G$
into the automorphisms of the star product. 
 We come back to this notion of covariance
in \S \ref{genmommap}.

A program to study  representations through deformation
quantization methods   was introduced  in 1978 by Bayen, 
Flato, Fronsdal, Lichnerowicz
and Sternheimer in \cite{BayenB, Stern1} and 
in Fronsdal \cite{Fronsdal}. The notion of covariant star 
products is essential for this.
The program involves a beautiful notion of adapted 
Fourier transform (or $*$ exponential).
I shall not speak about this aspect of 
deformation quantization; let me just mention that 
the program to study  representations through deformations 
was carried  first in the case
of nilpotent Lie groups in the mid eighties 
by Arnal and Cortet  \cite{ArnCor} and  that the 
adapted Fourier transform they found \cite{ArnCor2} gives a Kirillov type
formula for the character of the representation and a localisation
of the (packets of) coefficients of representations 
on the corresponding orbits
in the dual of the Lie algebra.  
This theory of star representations was developed for all representations
of compact groups and exponential groups,
and for some series of representations of semi-simple groups.

\medskip
$\bullet$ The standard star product on $\g^*$ does not
restrict to orbits (except for the Heisenberg group) so 
other algebraic constructions of star products on $S(\g)$ were
considered (with Michel Cahen in \cite{CG}, with Cahen 
and Arnal in \cite{ACG},
by Arnal, Ludwig and  Masmoudi in \cite{AL}
 and more recently by Fioresi and Lledo in \cite{FL}).
For instance, when $\g$ is semisimple, if $\H$ is the space of harmonic 
polynomials and if $I_1,\dots I_r$ are generators of the space of invariant 
polynomials, then any polynomial $P\in S(\g)$ writes uniquely as a sum
$P=\sum_{a_1\dots a_r}I^{a_1}_1\dots I^{a_r}_rh_{a_1\dots a_r}$ where
$h_{a_1\dots a_r}\in \H$. One considers the isomorphism $\sigma'$ 
between $S(\g)$
and $U(\g)$ induced by this decomposition
$$
\sigma'(P)=\sum_{a_1\dots a_r}(\sigma(I_1)\circ)^{a_1}\dots 
(\sigma(I_r)\circ)^{a_r}\circ \sigma(h_{a_1\dots a_r}).
$$
This gives a star product on $S(\g)$ which is not defined by differential
operators.
In fact, with Cahen and Rawnsley, we proved \cite{CGR} that if 
$\g$ is semisimple,
there is no differential star product on any 
neighbourhood of $0$ in $\g^*$ such that
$C *u=C u$ for the quadratic invariant polynomial $C\in S(\g)$ and 
$\forall u
 \in S(\g)$ (thus no differential star product which is tangential
  to the orbits).

\medskip
$\bullet$ In 1983, De Wilde and Lecomte proved \cite {DL} that on 
any symplectic manifold
 there exists a differential star product. This was obtained by imagining
 a very clever generalisation of a homogeneity condition in the 
 form of building
 at the same time the star product and a special derivation of it. 
 A very nice 
 presentation of this proof appears in \cite{DeWilde}.
 Their technique works to prove the existence of a differential star product
 on a regular Poisson manifold \cite{Masmoudi}.

\medskip
$\bullet$ In 1985, but appearing only in the West 
in the nineties \cite{Fed},
 Fedosov gave a  recursive construction of a star product 
 on a symplectic manifold
 $(M,\omega)$ by 
     
     -considering the Weyl  bundle $\W$ which is the bundle on $M$
      in associative Weyl algebras $W$ associated to the principal bundle
      of symplectic frames;
     
     -building, from a symplectic connection $\nabla$ on $M$, 
     a covariant derivative
      on $\Gamma(\W)$, $\partial=\partial^\nabla+[r,.]~ 
      ( r\in\Gamma(T^*M\otimes\W))$, such that
       $\partial\circ\partial=0$;
       
     -identifying $\CiM[[\nu]]$ with 
     $\{ s\in \Gamma(\W)\mid \partial s=0\}$ and transferring
     the associative pointwize product of sections to an 
     associative product on $\CiM[[\nu]]$.
 
 In 1994, he extended Êthis result to
  give a recursive construction in the  context of regular Poisson manifold \cite{Fed2}.
  
  \medskip
$\bullet$ Independently, also using the framework of Weyl bundles,
 Omori, Maeda and Yoshioka \cite{OMY1} gave  an 
alternative proof of existence of a differential star product
 on a symplectic manifold, gluing local
 Moyal star products.

\medskip
$\bullet$  In 1986, Drinfeld \cite{Drinfeld}  proposed a program of
quantization of  Poisson-Lie groups.

One way to consider the problem is to study
deformations of the Hopf algebra $U(\g)$ of a finite dimensional Lie
 bialgebra $(\g,p)$, where one deforms the product and the coproduct,
 the deformation of the coproduct being driven by the cocycle $p$.
  By duality this is the quantization of formal Poisson-Lie groups.
 For the standard structures on semisimple Lie algebras, this was solved
 by Drinfeld; it is the construction of  classical {\bf quantum groups}. 
 Drinfeld showed that the quantization is preferred, i.e. the product in
 $U(\g)$ is unchanged (or, by duality, the coproduct on formal series
-- which are the functions on the formal Poisson-Lie group--  is unchanged).
Etingof and Kazhdan \cite{EtKaz} proved that one can deform $U(\g)$ for any  Lie
 bialgebra $(\g,p)$, so one can always quantize formal Poisson-Lie groups.
 They proved that the deformation is differential. Pinczon \cite{Pinczon}
 proved that any preferred quantization  of a formal Poisson-Lie group
 is differential.
 
\smallskip
Quantum groups 
became very popular and had fundamental applications. It is  not my goal
to describe those developments but they brought new interest in the question
of deformation quantization of general Poisson manifolds. 

Indeed, another way to consider the problem of quantization of
a Poisson-Lie group $(G,P)$ is to study deformations  
of the Hopf algebra $C^\infty(G)$,
where the deformation of the product of functions is driven by $P$.

In the papers  \cite{BFGP} by Bonneau, Flato, Gerstenhaber and Pinczon
and \cite{BP2} by  Bidegain and Pinczon, the duality between those two approaches
is proven in the framework of topological deformations. In particular, any
classical quantum group gives, by topological duality, a differential
deformation of $C^\infty(G)$.
In \cite{BP}, Bidegain and Pinczon have built a differential preferred
deformation of $C^\infty(G)$ for any Poisson Lie group $(G,P)$ with $G$
semisimple.
They also show that any preferred deformation of $C^\infty(G)$
is automatically differential.
Then  Etingof and Kazhdan \cite{EtKaz2} proved that there exists
a differential deformation of $C^\infty(G)$ for any Poisson Lie group $(G,P)$,
deformation which is preferred in the quasi triangular case.

\medskip
$\bullet$  
The existence of a star product on other classes of Poisson manifolds
was studied by various authors (Omori--Maeda--Yoshioka \cite{OMY2}, 
Tamarkin \cite{Tam1}, Asin \cite{Asin},...).
 
\medskip
$\bullet$  In 1997, Kontsevich \cite{K} gave a proof 
of the existence of a star product on any
  Poisson manifold and gave an explicit formula 
  for a star product for any Poisson structure
  on $V=\R^m$. This appeared as a consequence of the proof 
 of his formality theorem. Tamarkin \cite{Tam2}  gave a
 version of the proof in the framework of the theory of operads.

  \subsection{Kontsevich star product on  $V$} 

Let $M$ be a domain in $V=\R^m$ and $\PO$ be any Poisson structure on $M$.
Kontsevich builds a star product on $(M,\PO)$, where the star product
 of two functions $u$ and $v$ 
is given in terms of some universal polydifferential operators
applied to the coefficients of the bi-vector field $\PO$ and 
to the functions $u,v$; the formula is invariant under
affine transformations of $\R^m$ and
the description of the $k^{\rm th}$ cochain
uses a special class $G_k$ of oriented labelled graphs.
      
An (oriented) graph $\Gamma$ is a pair $(V_{\Gamma},E_\Gamma)$ 
of two finite sets such that 
$E_\Gamma$ is a subset of $V_\Gamma\times V_\Gamma$;
elements of $V_\Gamma$ are  vertices of $\Gamma$, elements
of $E_\Gamma$ are edges of $\Gamma$. If 
$e=(v_1,v_2)\in E_\Gamma\subseteq V_\Gamma
\times V_\Gamma$ is an edge it is said that 
$e$ starts at $v_1$ and ends at $v_2$.
         
A graph $\Gamma$ belongs to $G_k$ if
$\Gamma$ has $k+2$ vertices $V_\Gamma=\{1,\dots,k\}\sqcup \{L,R\}$
and $2k$ labelled edges (with no multiple edges and no edge of the form
$(v,v)$ for $v \in V_\Gamma$),
$E_\Gamma=\{e_1^1,e_1^2,e_2^1,e_2^2,\dots,e_k^1,e_k^2\}$
where $e_j^1$ and $e_j^2$ start at $j$. 
       
To each labelled graph $\Gamma\in G_k$ and each skew symmetric $2$-tensor $P$
is associated a bidifferential operator $C_{\Gamma}(P)$  on $M\subset\R^m$ :
\begin{eqnarray*}
(C_{\Gamma}(P))(u,v):&=&\sum_{
I:E_\Gamma\longrightarrow \{1,\dots,m\}}[
 \prod_{j=1}^k (\prod_{e\in E_\Gamma\,\mid \, e=(*,j)} \partial_{I(e)})~~
 P^{I(e_j^1)I(e_j^2)}]\\
&\times& (\prod_{ e\in E_\Gamma\,\mid\, e=(*,L)} \partial_{I(e)})~~u\times
  (\prod_{e\in E_\Gamma\,\mid\, e=(*,R)} \partial_{I(e)})~~v.
\end{eqnarray*}
       
A weight $w_{\Gamma}\in \R$ is associated with each graph $\Gamma\in G_k$.
Denote by $\H_k$ the  space of configurations of $k$ numbered 
distinct points in the upper half plane $\H=\{z\in \C|\,\,Im(z)>0\}$:
$$
\H_k=\{(p_1,\dots,p_k)\mid p_j\in \H \quad 
           p_i\ne p_j\quad{\rm for} \quad i\ne j\}.
$$
For $p\in\H~,\, q\in \H\sqcup \R$ define
$$
\phi^{h}(p,q)={1\over 2 i} Log\left(\frac{(q-p)({\overline q}-p)}
              {(q-{\overline p})({\overline q} -
              {\overline p})}\right).
$$
 If $\Gamma\in G_k$ is a graph as above, assign a point
$p_j\in \H$ to the vertex  $j$ for 
$1\le j\le k$, point $0\in \R\subset \C$ to the vertex 
$L$, and point $1\in \R\subset \C$ to the vertex $R$. 
Every  edge $e\in E_\Gamma$ defines an ordered pair $(p,q)$ 
of points on $\H\sqcup \R$, thus a function
$\phi^{h}_e=\phi^h(p,q)$ on $\H_k$ with values in $\R/2 \pi \Z$.
              
The weight  of $\Gamma$ is defined by 
$$
w_\Gamma=\frac{1}{ k! (2\pi)^{2k}}\int_{\H_n}\wedge_{i=1}^k 
(d\phi^{h}_{e_i^1}\wedge d\phi^{h}_{e_i^2}).
$$

\begin{thm}(Kontsevich \cite{K})
Let $P$ be a Poisson tensor on a domain of $\R^m$. Then 
\[
u*_\nu v :=\sum_{k=0}^{\infty} {\nu^k}
    \sum_{\Gamma\in G_k} w_\Gamma (C_\Gamma (P))(u,v)
\]
defines a differential star product.
\end{thm}

\begin{definition} A star product on $\R^m$ is  {\bf universal}
if it is given, for any Poisson tensor $P$ by
\[
u*_\nu v=uv+\sum_{n>0}\nu^n (C_n(P))(u,v)\quad \rm{with} \quad C_n(P)=
\sum_{\Gamma \in G}w'(\Gamma) C_\Gamma(P)
\]
where the $w'(\Gamma)$ are scalars which are independent of $P$.
\end{definition}

\begin{prop}
Two universal star products $*_\nu$ and $*'_\nu$ on $\g^*$ are always 
equivalent modulo
a change of parameter. This means that there exists
a series $T(\nu) = \Id + \sum_{r = 1}^\infty \nu^r T_r$
of linear universal differential operators and a series
$\nu'= \nu +\sum_{n>1}f_n\nu^n$ such that
$$
u*'_{\nu'}v= T^{-1}(\nu)(T(\nu)u*_\nu T(\nu)v).
$$
\end{prop}

This was first proven by Arnal \cite{Arnal};
 we give here an elementary proof of this fact.
\begin{proof}
Assume by induction that $*$ and $*'$ coincide at order $n-1$.
Associativity relation at order $n$ implies 
$\partial C_n(P)=\partial C'_n(P)$; thus
$C_n(P)=C'_n(P)+\partial E_n(P) + A_n(P)$
where $E_n$ is a universal differential operator and $A_n(P)$
is universal, skewsymmetric and of order $(1,1)$.
Clearly,  if $A_n$ corresponds to a graph in $G_k$, 
it has $2k$ arrows and $k+2$
vertices, but since $P$ is linear on $\g^*$, at most one
arrow can end at any of the first $k$ vertices and since $A_n$
is $1$-differential in each argument, exactly one arrow ends
at each of the last $2$ vertices. Hence $2k\le k+2$ so $k\le 2$;
but if $k=1$ the only graph yields the Poisson bracket of functions
so this can be cancelled by a change of parameter and if
$k=2$ the graph corresponds to a symmetric bidifferential operator,
 hence does not
yield a $A_n(P)$. So we can assume that $A_n$ vanishes but then
the equivalence through $T(\nu)= \Id + \nu^n E_n(P)$ builds a star
product which is universal and coincide with $*$ at order $n$. 
\end{proof}

\begin{remark}
When a covariance condition  or a homogeneity
condition ($C_n(sP)=s^nC_n(P)$) is added, clearly there 
is no need for a change
 of parameter since the graph in $G_1$ can only arise in 
 $C_1$ so two such universal
 star products are always equivalent.

 In particular, the star product of Kontsevich and the standard 
 (CBH) star product
on $\g^*$, which in general are not the same, are equivalent.  
 The equivalence  is given by universal differential operators,
  hence combinations
of wheels which are graphs consisting of $k+1$ points 
$1,\dots,k$ and $L$ and $2k$
arrows $(1,2),(2,3)\dots  ,(k-1,k),(k,1)$ and $(j,L)\,1\le j\le k$.
 Such a wheel clearly vanishes when the Poisson structure corresponds
to a nilpotent Lie algebra so there is only one covariant universal 
star product on the
dual of a nilpotent Lie algebra (this appears in 
Arnal \cite{Arnal1} and Kathotia \cite{Kathotia}). 

The equivalence between Kontsevich and CBH star product has been explicitly
 constructed (see Arnal \cite{Arnal} and Dito \cite{Dito}) and gives
an integral formula for Kontsevich star product:
\[
u*v(\xi)=\int \tilde{u}(X)\tilde{v}(Y)\frac{F(X)F(Y)}{F(CBH(X,Y))}
   \exp(2\pi i \langle \xi,CBH(X,Y)\rangle) dXdY
   \]
where $F$ is some formal function on $\g$ written as a sum
of products of traces of powers of $\ad(X)$. 
   
 This star product has been used recently by Andler, 
 Dvorsky and Sahi \cite{Andler} to 
 establish a conjecture of Kashiwara and Vergne which, 
 in turn, gives a new proof of
 Duflo's result on the local solvability of bi-invariant 
 differential operators
 on a Lie group.
\end{remark}

\section{Equivalence of star products}

\begin{definition}
Two star products $*$ and $*'$ on $(M,P)$ are said to
be \textbf{equivalent} if there is a series
\[
T = \Id + \sum_{r = 1}^\infty \nu^r T_r
\]
where the $T_r$ are linear operators on $\CiM$, such that
\begin{equation}\label{intro:equivalence}
T(f*g) = Tf*' Tg.
\end{equation}
\end{definition}

\noindent Remark that the $T_r$ automatically vanish on
constants since $1$ is a unit for $*$ and for $*'$. Using in a similar way 
linear operators which do not necessarily vanish on constants, one can
pass from any  associative deformation of the product of functions
on a Poisson manifold $(M,P)$ to another such deformation with $1$ being a unit.  

\begin{definition}
A {\textbf{Poisson deformation}} of the Poisson bracket on a Poisson manifold
$(M,P)$ is a Lie algebra deformation of  $(\CiM,\{~,~\})$ which is a derivation
 in each argument, i.e. of the form $\{u,v\}_\nu=P_\nu(du,dv)$
 where $P_\nu=P+\sum\nu^kP_k$ is a series of skewsymmetric 
 contravariant $2$-tensors on $M$ (such that $[P_\nu,P_\nu]=0$).
 
Two Poisson deformations  
$P_\nu$ and $P'_\nu$ of the Poisson bracket $P$
on a Poisson manifold $(M,P)$  are {\textbf{equivalent}}
if there exists  a formal path in the diffeomorphism group of $M$, 
starting at the identity, i. e. a series 
$T=\exp D=\Id+\sum_j \frac{1}{j!} D^j$
for $D=\sum_{r\ge 1}\nu^r D_r$ where the 
$D_r$ are  vector fields on $M$, such that 
$$
T\{u,v\}_\nu=\{Tu,Tv\}'_\nu
$$
where  $\{u,v\}_\nu=P_\nu(du,dv)$ and $\{u,v\}'_\nu=P'_\nu(du,dv)$.
\end{definition}

\medskip
 $\bullet$ In the general theory of deformations, Gerstenhaber
  \cite{Gerst} showed
  how equivalence is linked to some second cohomology space.

\medskip
 $\bullet$ For symplectic manifolds,  Flato, Lichnerowicz 
 and Sternheimer in 1974
studied $1$-differential 
deformations of the Poisson bracket \cite{FLS}; it follows from their work,
and appears  in Lecomte \cite{Lecomte}, that:
\begin{prop}
The equivalence classes of Poisson deformations of the Poisson bracket $P$
on a symplectic manifold $(M,\omega)$ are  parametrised by 
 $H^2(M;\R)[[\nu]]$.
\end{prop}
Indeed, one first show that
any Poisson deformation $P_\nu$ of the Poisson bracket $P$
on a symplectic manifold $(M,\omega)$ is of the form $P^\Omega$
for a series $\Omega=\omega+\sum_{k\ge 1}\nu^k\omega_k$
where the $\omega_k$ are closed $2$-forms, and
$P^\Omega(du,dv)=-\Omega(X^\Omega_u,X^\Omega_v)$
where $X^\Omega_u=X_u+\nu(\dots) \in \Gamma (TM) [[\nu]]$ is
 the element defined by
$i(X^\Omega_u)\Omega=du$. 

\begin{itemize}
\item[~]Observe that for any series  $\Omega$
as above the series of $1$-differential $2$-cochains 
$P^\Omega$ satisfies $[P^\Omega,P^\Omega]=0$ because $\Omega$ 
is a closed $2$-form,
so defines indeed a Poisson  deformation of $P$.  

Reciprocally, given a 
$1$-differential deformation $P_\nu=P+\sum_j \nu^jP_j$,  assume
 it coincides up to order $k$ ($k\ge 0$) with $P^{\Omega_k}$ for some
$\Omega_k=\omega+\sum_{j\ge 1}^k\nu^j\omega_j$ then
$[P_\nu,P_\nu]-[P^{\Omega_k},P^{\Omega_k}]=0$ imply at order
$k+1$ that $[P,P_{k+1}-(P^{\Omega_k})_{k+1}=0]$ so that
there exists a closed $2$-form $\omega_{k+1}$
such that $(P_{k+1}-(P^{\Omega_k})_{k+1})(du,dv)=\omega_{k+1}(X_u,X_v)$.
Hence $P_\nu$ coincides up to order $k+1$ with $P^{\Omega_{k+1}}$ where
$\Omega_{k+1}=\omega+\sum_{j\ge 1}^{k+1}\nu^j\omega_j$. By induction,
any Poisson deformation of $P$ is of the form $P^\Omega$ for
a series $\Omega=\omega+\sum_{k\ge 1}\nu^k\omega_k$
where the $\omega_k$ are closed $2$-forms.
\end{itemize}
One then shows that two  Poisson deformations $P^\Omega$ and $P^{\Omega'}$
are equivalent if and only if $\omega_k$ and $\omega'_k$ are cohomologous
for all $k\ge 1$.

\begin{itemize}
\item[~]
If $\Omega$ and $\Omega'$ coincide to order $k-1$ and  
$\omega_k-\omega'_k=dF_k$,
then  $\{u,v\}_\nu$ and $T^{-1}\{ Tu,Tv\}'_\nu$, with 
$T= \exp D\quad Du=\nu^{k-1} F_k(X_u)$,
correspond to forms which coincide to order $k$ 
so one has equivalence
if all forms are cohomologous. 

Reciprocally, if $P^\Omega$ and $P^{\Omega'}$ are equivalent and
$\Omega$ and $\Omega'$ coincide to order $k-1$, then 
the equivalence can be written
$T= \exp D\quad Du=\sum_{j\ge k}\nu^{j} F_j(X_u)$ and 
$\omega_k-\omega'_k=dF_k$.
Indeed if $Tu=u +\nu^r F_r(X_u) +\dots$ , then the relation
$P^\Omega(du,dv)=T^{-1}P^{\Omega'}(dTu,dTv)$ at order 
$r\le k$ yields $\omega_r-\omega'_r=dF_r$.
When $r<k$ this means $dF_r=0$; locally, on a contractible set  
$U$, $F_{r\vert U}=df_U$. The map 
$D_r: u\mapsto P^{\Omega'}(df_U,du)$ is globally defined and is a derivation
of $\{~,~\}'$ so  $T'=\exp -\nu^r D_r\circ T=\Id + \nu^{r+1}F_{r+1}+\dots$
is still an equivalence. By induction, one gets the result.
\end{itemize}

\medskip
 $\bullet$In 1978, Bayen, Flato, Fronsdal, Lichnerowicz and 
 Sternheimer \cite{BayenB2}
stressed that different
orderings in physics lead to equivalent star products on $\R^{2n}$.
This shows that the notion of mathematical equivalence is not the same as the
notion of physical equivalence (i.e. two star products 
leading to the same spectrum for 
each observable); we studied this
difference  with Cahen, Flato and Sternheimer in \cite{CFGS}. 
Bayen et al.~ also proved in \cite{BayenB} that Moyal 
star product is the only star product
whose cochains are given by  polynomials
in the Poisson structure $P$; this was the first 
consideration of some universality
property to build and classify star products.

\medskip
$\bullet$ Recall that a star product $*$ on $(M,\omega)$ 
is called differential if the
$2$-cochains $C_r(u,v)$ giving it are bi-differential operators.
As was observed by Lichnerowicz \cite{Lichne} 
and  Deligne \cite{Deligne}~: 
\begin{prop}
If $*$ and $*'$ are differential star products and $T(u) = u +
\sum_{r\ge1} \nu^r T_r(u)$ is an equivalence so that $T(u*v) =
T(u)*'T(v)$ then the $T_r$ are differential operators.
\end{prop}
\begin{proof}
Indeed if $T=\Id + \nu^kT_k+\dots$ then $\cob T_k=C'_k-C_k$
is differential so $C'_k-C_k$ is a differential $2$-cocycle
with vanishing skewsymmetric part but then, using Vey's formula,
 it is the coboundary of a differential $1$-cochain $E$ and
 $T_k-E$, being a $1$-cocycle, is a vector field so $T_k$ is differential.
 One then proceeds by induction, considering 
 $T'= (\Id + \nu^kT_k)^{-1}\circ T=\Id +\nu^{k+1}T'_{k+1}+\dots$ and
 the two  differential star products $*$ and $*''$, 
 where $u*''v=(\Id + \nu^kT_k)^{-1}((\Id + \nu^kT_k)u*'(\Id + \nu^kT_k)v),$ 
  which are equivalent through $T'$ (i.e. $T'(u*v) =T'(u)*''T'(v)$).
  \end{proof} 

\medskip
$\bullet$ A differential star product is equivalent to one 
with linear term in
$\nu$ given by $\frac12 \{u,v\}$. Indeed $C_1(u,v)$
is a Hochschild cocycle with antisymmetric part given by $\frac12
\{u,v\}$ so $C_1 = \frac12 P + \cob B$ for a differential $1$-cochain $B$.
Setting $T(u) = u +
\nu B(u)$ and $u*'v = T(T^{-1}(u)*T^{-1}(v))$,  this equivalent star product
$*'$ has the required form.

\medskip
 $\bullet$ In 1979, we proved \cite{Guequiv} that all differential
 deformed brackets on  $\R^{2n}$  (or on any
symplectic manifold such that $b_2=0$) are equivalent 
modulo a change of the parameter,
and this implies a similar result for star products; this was proven
by direct methods by Lichnerowicz \cite{Lichneequiv}:
\begin{prop}\label{equib2}
Let $*$ and $*'$ be two differential star products on $(M,\omega)$ and
suppose that $\coh2=0$. Then there exists a local equivalence $T = \Id +
\sum_{k\ge1}\nu^kT_k$ on $\CiM\bbnu$ such that $u*'v = T(T^{-1}u*T^{-1}v)$
for all $u, v \in \CiM\bbnu$.
\end{prop}

\begin{proof}
Let us suppose that, modulo some equivalence, the two star products $*$
and $*'$ coincide up to order $k$. Then associativity at order $k$ shows
that $C_k-C'_k$ is a Hochschild $2$-cocycle and so by
(\ref{cohomoVey}) can be written as $(C_k-C'_k)(u,v)=(\cob B)(u,v) +
A(X_u,X_v)$ for a $2$-form $A$. The total skewsymmetrization 
of the associativity
relation at order $k + 1$ shows that $A$ is a closed $2$-form. Since the
second cohomology vanishes, $A$ is exact, $A = dF$. Transforming by the
equivalence defined by $Tu = u + \nu^{k-1}2F(X_u)$, we can assume that
the skewsymmetric part of $C_k-C'_k$ vanishes. Then
$C_k-C'_k = \cob B $ where $B$ is a
differential operator. Using the equivalence defined by $T =
I + \nu^k B$ we can assume that the star products coincide, modulo an
equivalence, up to order $k + 1$ and the result follows from induction
 since two star products always agree in their leading term.
\end{proof}

\medskip
$\bullet$ It followed from 
the above proof and  results similar  to \cite{Guequiv} 
(i.e. two star products which are equivalent
 and coincide at order $k$ differ at order $k+1$ by a Hochschild 
 2-cocycle whose skewsymmetric part corresponds to an exact 2-form)
 that at each step in $\nu$, equivalence
classes of differential star products on a symplectic manifold
$(M,\omega)$ are parametrised by 
$H^2(M;\R)$, if all such deformations exist. 
The general existence was proven by De Wilde and Lecomte. 
At that time, one assumed the
parity condition $C_n(u,v)=(-1)^nC_n(v,u)$, so equivalence 
classes of such differential
star products were parametrised by series $H^2(M;\R)[[\nu^2]]$. 
The parametrization
was  not  canonical.

\medskip
$\bullet$ In 1994, Fedosov proved that his recursive 
construction works in a 
more general setting : given any series of closed 
2-forms on a symplectic manifold $(M,\omega)$,
he could build a connection on the Weyl bundle whose curvature is linked
to that series and   a star product whose equivalence
class only depends on the element in $H^2(M;\R)[[\nu]]$ 
corresponding to that
 series of forms.
 
 \medskip
 $\bullet$ In 1995, Nest and Tsygan \cite{Nest-Tsygan},
  then Deligne \cite{Deligne} and Bertelson 
 \cite{Bertelson}   proved 
 that any differential star product on a symplectic manifold $(M,\omega)$
 is equivalent to a Fedosov star product and that its  equivalence class
is parametrised by the corresponding element in $H^2(M;\R)[[\nu]]$.

\medskip
$\bullet$ In 1997, Kontsevich \cite{K} proved that the coincidence
of the set of  equivalence classes of star and Poisson deformations 
is true for general Poisson manifolds :
\begin{thm}
The  set of  equivalence classes of differential star products on a
 Poisson manifold $(M,\PO)$ can be naturally identified with
the set of equivalence classes of Poisson deformations of $\PO$:
 $$\PO_{\nu}=\PO\nu+\PO_2 \nu^2+\dots\in \Gamma(X,\wedge^2 T_X)[[\nu]],\,\,\,
 [\PO_{\nu},\PO_{\nu}]=0 .$$  
\end{thm}

\medskip
$\bullet$ Remark that all results concerning parametrisation of 
equivalence classes of 
differential star products are still valid for star products 
defined by local cochains or for star products defined by continuous
cochains (\cite{Gutt},
 Pinczon \cite{Pinczon}).
Parametrization of  equivalence classes
 of special star products have been obtained :
  star products with separation of variables
(by Karabegov \cite{Karaequivsepvar}), 
invariant star products on a symplectic manifold
when there exists an invariant symplectic connection 
(with Bertelson and Bieliavsky
\cite{BBG}), algebraic star products 
 (Chloup \cite{Chloup}, Kontsevich \cite{K})...
 
 \medskip
 $\bullet$ The association of an element in $H^2(M;\R)[[\nu]]$
to  the equivalence class of a star product 
on a symplectic manifold is one way
to associate an invariant to a star product; 
other such associations are obtained
by star version of index theorems and trace functionals on the algebra 
$(\CiM[[\nu]],*)$. I shall not develop that aspect here. It was first
considered by Connes, Flato and Sternheimer in \cite{CFS} where they
introduce the notion of closed star product, i.e. such that
$$
\int_M a*b~ \omega^n=\int_M b*a ~\omega^n 
~~({\rm{mod }}~\nu^n) ~~\forall a,b\in\CiM[[\nu]], 
$$
and show how their classification is linked to cyclic cohomology.

They obtain, for the cotangent bundle to a compact Riemannian manifold,
its Todd class as the ``character'' associated to the star product 
corresponding to normal ordering.

The notion of a trace for star products  
and the star version of index theorems
have been studied by Fedosov \cite{Fed2,Fedosov,Fed3}
 and by Nest and Tsygan \cite{Nest-Tsygan,NT2}. 
        
\subsection{Deligne's cohomology classes associated to differential
star products on symplectic manifolds}

 Deligne defines two cohomological classes associated to
differential star products on a symplectic manifold. 
This leads to an intrinsic
way to parametrise the equivalence class of such a differential star product.
Although the question makes sense more generally for Poisson manifolds,
Deligne's method depends crucially on the Darboux theorem and the
uniqueness of the Moyal star product on ${\R}^{2n}$ so the methods do
not extend to general Poisson manifolds. 

The first class is
a relative class; fixing a star product on the manifold, it
intrinsically associates to any equivalence class of star products an
element in $\coh2\bbnu$. This is done in \cech\ cohomology by
looking at the obstruction to gluing local equivalences.

Deligne's second class is built from special local derivations of a star
product. The same derivations played a special role in the first general
existence theorem \cite{DL} for a star product on a symplectic
manifold. Deligne used some properties of Fedosov's construction and
central curvature class to relate his two classes and to see how to
characterise an equivalence class of star products by the derivation
related class and some extra data obtained from the second term in the
deformation. With John Rawnsley \cite{GR}, we did
 this by direct \cech\ methods which I shall present here.

\subsubsection{The relative class}

Let $*$ and $*'$ be two differential star products on $(M,\omega)$.
Let $U$ be a contractible open subset of $M$ and
$N_U = C^\infty(U)$. Remark that any differential star product on $M$
restricts to $U$ and $\coh2(U)=0$, hence, by proposition \ref{equib2},
there exists a local equivalence $T= \Id +
\sum_{k\ge1}\nu^kT_k$ on $N_U\bbnu$
so that $u*'v = T({T}^{-1}u*{T}^{-1}v)$ for all
$u, v \in N_U\bbnu$.

\begin{prop}
Let $*$ be a differential star product on $(M,\omega)$ and
suppose that $\coh1$ vanishes.
\begin{itemize}
\item
Any self-equivalence $A = \Id + \sum_{k\ge1} \nu^k A_k$ of
$*$ is inner: $A =
\exp\, \ad_* a$ for some $a \in \CiM\bbnu$.
\item
Any $\nu$-linear derivation of $*$  is of the form
$ D=\sum_{i\ge 0}\nu^i D_i $
where each $D_i$ corresponds to a symplectic vector field $X_i$
and is given on a contractible open set $U$ by
$$
D_i u|_U=\frac{1}{ \nu} \left(f_i^U*u-u*f_i^U\right)
$$
if $X_iu|_U = \{ f_i^U,u\}|_U$.
\end{itemize}
\end{prop}
Indeed, one builds $a$ recursively; 
assuming $A= \Id + \sum_{r\ge k} \nu^r A_r$ and $k\ge 1$,
the condition $A(u*v) = Au*Av$ implies at order $k$ in $\nu$
 that $A_k(uv) + C_k(u,v) = A_k(u)v
+ uA_k(v) + C_k(u,v) $ so that $A_k$ is a vector field. Taking the
skew part of the terms in $\nu^{k+1}$ we have that $A_k$ is a derivation
of the Poisson bracket. Since $\coh1=0$, one can write 
$A_k(u) = \{a_{k-1},u\}$ for some function $a_{k-1}$. 
 Then $(\exp\,-\ad_* \nu^{k-1}a_{k-1}) \circ A = \Id +
O(\nu^{k+1})$ and the induction proceeds.
The proof for $\nu$-linear derivation is similar.

\medskip

The above results  can  be applied
to the restriction of a differential star product on $(M,\omega)$
to a contractible open set $U$. Set, as above, $N_U = C^\infty(U)$.
If $A = \Id + \sum_{k\ge1} \nu^k A_k$ is a formal linear
operator on $N_U\bbnu$ which preserves the differential
star product $*$, then there is $a \in N_U\bbnu$ with
$A = \exp\, \ad_* a$.
 Similarly, any local
$\nu$-linear derivation $D_U$ of   $*$ on $N_U\bbnu$ 
 is essentially inner: $D_U =
\frac1\nu \ad_* d_U$ for some $d_U \in N_U\bbnu$.

\medskip
It is convenient to write the composition of automorphisms of the
form $\exp\ad_* a$ in terms of $a$.
In a pronilpotent situation this is done with the
\textbf{Campbell--Baker--Hausdorff composition} which is denoted by
$a \circ_*b$:
$$
a \circ_*b=a + \int_0^1 \psi(\exp\,\ad_* a\circ \exp\,t\ad_* b)b\, dt
$$
where
$$
\psi(z)=\frac{z\log(z)}{z-1}=\sum_{n\ge 1}\left( \frac{(-1)^n}{n + 1}
+ {\frac{(-1)^{n + 1}}{n}}\right) (z-1)^n.
$$
Notice that the formula is well defined
(at any given order in $\nu$, only a finite number of terms
arise) and it is given by the usual series
\[
a \circ_* b = a + b + \frac12 [a,b]_* +
\frac{1}{12}( [a,[a,b]_*]_* + [b,[b,a]_*]_*) \cdots.
\]
The following results are standard (N.~Bourbaki, Groupes et alg\`ebres de Lie,
\textit{\'El\'ements de Math\'ematique}, Livre 9, Chapitre 2, \S6):
\begin{itemize}
\item $\circ_*$ is an associative composition law;
\item $\exp\,\ad_* (a \circ_* b) = \exp\,\ad_* a \circ \exp\,\ad_* b$;
\item $a \circ_* b \circ_*(-a) = \exp( \ad_* a)\, b$;
\item $-(a \circ_* b) = (-b) \circ_* (-a)$;
\item $\displaystyle{\left.\frac{d}{dt}\right|_0 (-a)\circ_*(a + tb) =
\frac{1 - \exp \,(-\ad_*a)}{\ad_*a} (b)}$.
\end{itemize}

Let $(M,\omega)$ be a symplectic manifold. We fix a locally finite open
cover $\mathcal{U} = \{U_\alpha\}_{\alpha \in I}$ by Darboux coordinate
charts such that the $U_\alpha$ and all their non-empty intersections
are contractible, and we fix a partition of unity
$\{\theta_\alpha\}_{\alpha \in I}$ subordinate to $\mathcal{U}$. Set
$N_\alpha = C^\infty(U_\alpha)$, $N_\alphabeta = C^\infty(U_\alpha \cap
U_\beta)$, and so on.

Now suppose that $*$ and $*'$ are two differential star products on
$(M,\omega)$. We have seen that their restrictions to $N_\alpha\bbnu$
are equivalent so there exist formal differential operators $ T_\alpha
\colon N_\alpha\bbnu \to N_\alpha\bbnu$ such that
\[
T_\alpha (u*v) = T_\alpha (u) *' T_\alpha (v), \qquad
u,v \in N_\alpha\bbnu.
\]
On $U_\alpha\cap U_\beta$, $T_\beta^{-1} \circ T_\alpha$ will be a
self-equivalence of $*$ on $N_\alphabeta\bbnu$ and so there will be
elements $t_\betaalpha = - t_\alphabeta$ in $N_\alphabeta\bbnu$ with
\[
T_\beta^{-1} \circ T_\alpha = \exp\,\ad_* t_\betaalpha.
\]
On $U_\alpha\cap U_\beta\cap U_\gamma$ the element
$$
t_\gammabetaalpha = t_\alphagamma \circ_* t_\gammabeta \circ_*
t_\betaalpha
$$
induces the identity automorphism and hence is in the
centre $\R\bbnu$ of $N_\alphabetagamma\bbnu$. The family of
$t_\gammabetaalpha$ is thus a \cech\ $2$-cocycle for the covering
$\mathcal{U}$ with values in $\R\bbnu$. The standard arguments show
that its class does not depend on the choices made, and is
compatible with refinements. Since every open cover has a refinement
of the kind considered it follows that $t_\gammabetaalpha$
determines a unique \cech\ cohomology class $[t_\gammabetaalpha]
\in\coh2\bbnu$.

\begin{definition}
\[
t(*',*) = [t_\gammabetaalpha] \in \coh2\bbnu
\]
is \textbf{Deligne's relative class}.
\end{definition}

It is easy to see, using the fact that the cohomology
of the sheaf of smooth functions is trivial:

\begin{thm}(Deligne)
Fixing a differential star product $*$, the class
$t(*',*)$ in $\coh2\bbnu$ depends only on the equivalence class of
the differential star product $*'$, and sets up a bijection
between the set of equivalence classes of differential star
products and $\coh2\bbnu$.

If $*$, $*'$, $*''$ are three differential star products on
$(M,\omega)$ then
\begin{equation}\label{relative:addition}
t(*'',*) = t(*'',*') + t(*',*).
\end{equation}

\end{thm}

\subsubsection{The  derivation related class}

The addition formula above suggests that
$t(*',*)$ should be a difference of classes $c(*'), c(*) \in
\coh2\bbnu$. Moreover, the class $c(*)$
should determine the star product $*$ up to equivalence. 
\begin{definition}
Let $U$ be an open set of $M$.
Say that a derivation $D$ of $(C^\infty(U)\bbnu,*)$ is \textbf{$\nu$-Euler} 
 if it has the form
\begin{equation}\label{intrinsic:dereqn}
D = \nu \frac{\partial}{\partial\nu} + X + D'
\end{equation}
where $X$ is conformally symplectic on $U$ 
($\L_X\omega\vert_U=\omega\vert_U$) and
$D' = \sum_{r\ge1} \nu^r D'_r$ with the $D'_r$ differential
operators on $U$.
\end{definition}

\begin{prop}\label{intrinsic:derivations}
Let $*$ be a differential star product on $(M,\omega)$. For each
$U_\alpha \in \mathcal{U}$ there exists a $\nu$-Euler derivation $D_\alpha =
\nu \frac{\partial}{\partial\nu} + X_\alpha + D'_\alpha$ of the algebra
$(N_\alpha\bbnu,*)$.
\end{prop}

\begin{proof}
On an open set in $\R^{2n}$ with the standard symplectic structure
$\Omega$, denote the Poisson bracket by $P$. Let $X$ be a conformal
vector field so $\L_X\Omega = \Omega$. 
The Moyal star product $*_M$ is given by
$u*_M v = uv + \sum_{r\ge1} ({\frac{\nu}{2}})^r/r! P^r(u,v)$ and
 $D = \nu \frac{\partial}{\partial\nu} + X$ is a derivation
of $*_M$.

Now $(U_\alpha, \omega)$ is symplectomorphic to an open set in $\R^{2n}$ and
any differential star product on this open set is equivalent to $*_M$ so
we can pull back $D$ and $*_M$ to $U_\alpha$ by a symplectomorphism
to give a star product $*'$ with a derivation of the form $\nu
\frac{\partial}{\partial\nu} + X_\alpha$. If $T$ is an equivalence of
$*$ with $*'$ on $U_\alpha$ then $D_\alpha = T^{-1} \circ (\nu
\frac{\partial}{\partial\nu} + X_\alpha) \circ T$ is a derivation of the
required form.
\end{proof}

We take such a collection of derivations $D_\alpha$ given by
Proposition \ref{intrinsic:derivations} and on $U_\alpha \cap
U_\beta$ we consider the differences $D_\beta - D_\alpha$. They are
derivations of $*$ and the $\nu$ derivatives cancel out, so $D_\beta
- D_\alpha$ is a $\nu$-linear derivation of $N_\alphabeta\bbnu$. Any
$\nu$-linear derivation is of the form $\frac1\nu \ad_* d$, so there
are $d_\betaalpha \in N_\alphabeta\bbnu$ with
\begin{equation}\label{intrinsic:difference}
D_\beta - D_\alpha = \frac1\nu \ad_* d_\betaalpha
\end{equation}
with $d_\betaalpha$ unique up to a central element. On $U_\alpha \cap
U_\beta \cap U_\gamma$ the combination $d_\alphagamma + d_\gammabeta +
d_\betaalpha$ must be central and hence defines $d_\gammabetaalpha \in
\R\bbnu$. It is easy to see that $ d_\gammabetaalpha$ is a $2$-cocycle
whose \cech\ class $[d_\gammabetaalpha] \in \coh2\bbnu$ does not depend
on any of the choices made.

\begin{definition}
$d(*) = [d_\gammabetaalpha] \in \coh2\bbnu$ is \textbf{Deligne's intrinsic
derivation-related class}.
\end{definition}

\begin{itemize}
\item
In fact the class considered by Deligne is actually $\frac1\nu d(*)$.
 A purely \cech-theoretic accounts of
this class is given in Karabegov \cite{Karaequivsepvar}.
\item
If $*$ and $*'$ are equivalent  then
$d(*') = d(*)$.
\item
If $d(*) = \sum_{r\ge0} \nu^r d^r(*)$ then $d^0(*) = [\omega]$ under
the de~Rham isomorphism, and $d^1(*) = 0$.
\end{itemize}
Consider two differential star products $*$ and $*'$ on $(M,\omega)$
with local equivalences $T_\alpha$ and local $\nu$-Euler derivations
$D_\alpha$ for $*$ . Then $D'_\alpha = T_\alpha \circ D_\alpha \circ
T_\alpha^{-1}$ are local $\nu$-Euler derivations for $*'$. Let
$D_\beta - D_\alpha = \frac1\nu \ad_* d_\betaalpha$ and $T_\beta^{-1}
\circ T_\alpha = \exp\,\ad_* t_\betaalpha$ on $U_\alpha \cap U_\beta$.
Then $D'_\beta - D'_\alpha = \frac1\nu {\ad}_{*'} d'_\betaalpha$ where
\[
d'_\betaalpha = T_\beta d_\betaalpha - \nu T_\beta \circ \left(
\frac{1-\exp\,(-\ad_* t_{\alpha\beta})}{\ad_* t_{\alpha\beta}}\right)
\circ D_\alpha t_{\alpha\beta}.
\]
In this situation
\[
d'_\gammabetaalpha = T_\alpha(d_\gammabetaalpha + \nu^2
\frac{\partial}{\partial\nu} t_\gammabetaalpha).
\]
This gives a direct proof of:

\begin{thm}(Deligne)\label{intrinsic:relat}
The relative class and the intrinsic derivation-related classes of two
differential star products $*$ and $*'$ are related by
\begin{equation}\label{eq:intrel}
\nu^2 \frac {\partial}{\partial\nu} t(*',*) = d(*') - d(*).
\end{equation}
\end{thm}

\subsubsection{The characteristic class}

The formula above  shows that the
information which is ``lost'' in $d(*')-d(*)$ corresponds to the zeroth
order term in $\nu$ of $t(*',*)$. 

\begin{remark}
In \cite{Gutt23,DGL} it was shown that any bidifferential
operator $C$, vanishing on constants, which is a $2$-cocycle for the
Chevalley cohomology of $(\CiM,\{ ~,~\} )$ with values in $\CiM$ associated to
the adjoint representation (i.e.\ such that $
S_{u,v,w} [\{ u, C(v,w)\} -C(\{ u,v \} , w)] =0 $
where $S_{u,v,w}$ denotes the sum over cyclic permutations of $u,v$ and
$w$) can be written as
$$
C(u,v)=aS^3_\Gamma (u,v) + A(X_u,X_v) + [ \{ u,Ev\} + \{ Eu,v\}
-E(\{ u,v\} )]
$$
where $a\in\R$, where $S^3_\Gamma$ is a bidifferential $2$-cocycle
introduced in \cite{BayenB}
(which vanishes on constants and is never a coboundary and whose symbol is
of order $3$ in each argument), where $A$ is a closed $2$-form on $M$ and
where $E$ is a differential operator vanishing on constants. Hence
$$
H^2_{{\rm{Chev,nc}}}(\CiM,\CiM)=\R \oplus \coh2
$$
and we define the ${\#}$ operator as the projection 
on the second factor relative
to this decomposition.
\end{remark}

\begin{prop}
Given two differential star products $*$ and $*'$, the  term of order zero
in Deligne's relative class $t(*',*)=\sum_{r\ge 0} \nu^r t^r(*',*)$ is
given by
$$
t^0(*',*)=-2({C'}^-_2)^{\#} + 2(C^-_2)^{\#}.
$$
\end{prop}
If $C_1=\frac{1}{2}\{\, ,\,\}$, then $C_2^-(u,v)=A(X_u,X_v)$ where 
$A$ is a closed $2$-form and
$(C^-_2)^{\#}=[A]$ so it ``is'' the skewsymmetric part of $C_2$.

It follows from what we did before that the association to a
differential star product of $(C^-_2)^{\#}$ and $d(*)$ completely
determines its equivalence class.

\begin{definition}
The \textbf{characteristic class}  of a differential star product
$*$ on $(M,\omega )$ is the element $c(*)$ of the affine space ${\frac{-[
\omega ]}{\nu}} + \coh2\bbnu$ defined by
\begin{eqnarray*}
c(*)^0&=&-2(C^-_2)^{\#}\\
{\frac{\partial}{\partial \nu}} c(*)(\nu)&=&{\frac{1}{\nu^2}} d(*)
\end{eqnarray*}
\end{definition}

\begin{thm}
The characteristic class has the following properties:
\begin{itemize}
\item The relative class is given  by
\begin{equation}
t(*',*)=c(*')-c(*)
\end{equation}
\item The map $C$ from equivalence classes of star products on $(M,\omega)$
to the affine
space ${\frac{-[ \omega ]}{\nu}} + \coh2\bbnu$ mapping $[ * ] $ to $c(*)$
is a bijection.
\item If $\psi \colon M\to M'$ is a diffeomorphism and if $*$ is a
star product on $(M,\omega )$ then $u*'v=(\psi^{-1})^*(\psi^*u * \psi^*
v)$ defines a star product denoted $*'=(\psi^{-1})^* *$ on
$(M',\omega')$ where $\omega'=(\psi^{-1})^* \omega$. The characteristic
class is natural relative to diffeomorphisms:
\begin{equation}
c((\psi^{-1})^* *)=(\psi^{-1})^* c(*).
\end{equation}
\item Consider a change of parameter
$f(\nu)=\sum_{r\ge 1} \nu^r f_r$ where $f_r \in \R$ and $f_1\ne 0$
and let $*'$ be the star product obtained from $*$ by this change of
parameter, i.e.
$ u*'v=u.v + \sum_{r\ge 1} (f(\nu))^r C_r(u,v)=u.v + f_1\nu C_1(u,v)
+ \nu ^2 ((f_1)^2C_2 (u,v) + f_2 C_1(u,v)) + \ldots $. Then $*'$
is a differential star product on $(M,\omega')$ where $\omega'=
{\frac{1}{f_1}}\omega$ and we have equivariance under a change of
parameter:
\begin{equation}
c(*')(\nu)=c(*)(f(\nu)).
\end{equation}
\end{itemize}
\end{thm}

The characteristic class $c(*)$ coincides (cf Deligne \cite{Deligne}
and Neumaier \cite{Neumaier}) for Fedosov-type star products with their
characteristic class introduced by Fedosov as the de~Rham class of the
curvature of the generalised connection used to build them 
(up to a sign and factors of $2$). That characteristic class
is also studied by Weinstein and Xu in \cite{WX}.
 The fact that $d(*)$ and
$(C^-_2)^{\#}$ completely characterise the equivalence class of a star
product is also proven by \cech\ methods in De~Wilde
\cite{DeWilde}.

\subsection{Automorphisms of a star product 
and generalised moment map}\label{genmommap}

The above proposition allows to study  automorphisms
of star products  on a symplectic manifold (\cite{Rauch}, \cite{GR}).
\begin{definition}
An \textbf{isomorphism} from a differential star product $*$ on
$(M,\omega)$ to a differential star product $*'$ on $(M',\omega')$ is an
$\R$-linear bijective map $A \colon \CiM\bbnu \to C^\infty(M')\bbnu$, continuous in
the $\nu$-adic topology (i.e. $A(\sum_r \nu^r u_r)$ is the limit of
$\sum_{r\le N}A(\nu^r u_r)$ ), such that
$$
A(u*v)=Au*'Av.
$$
\end{definition}

Notice that if $A$ is such an isomorphism, then $A(\nu)$ is central for
$*'$ so that $A(\nu)=f(\nu)$ where $f(\nu)\in \R\bbnu$ is without
constant term to get the $\nu$-adic continuity. Let us denote by $*''$
the differential star product on $(M,\omega_1 = {\frac{1}{f_1}}\omega)$
obtained by a change of parameter
$$
u*^{''}_{\nu}v=u*_{f(\nu)}v=F(F^{-1}u*F^{-1}v)
$$
for $F \colon \CiM\bbnu \to \CiM\bbnu \colon \sum_r\nu^ru_r\mapsto \sum_r
f(\nu)^r u_r$. 

\noindent Define $A': \CiM\bbnu \to C^\infty(M')\bbnu$ by $A=A'\circ F$.
 Then $A'$ is a $\nu$-linear isomorphism between $*''$ and $*'$:
$$
A'(u*''v)=A'u*'A'v.
$$
At order zero in $\nu$ this yields $A'_0(u.v)=A'_0u.A'_0v$ 
so that there exists a diffeomorphism $\psi:M' \to M$ with $A'_0u=\psi^*u$.
The skewsymmetric part of the isomorphism relation at order $1$ in $\nu$
implies that $\psi^*\omega_1=\omega'$.
Let us denote by $*'''$ the differential star product on $(M,\omega_1)$
obtained by pullback via $\psi$ of $*'$:
$$
u*'''v=(\psi^{-1})^* (\psi^*u *' \psi^*v)
$$
and define $B: \CiM \bbnu \to \CiM \bbnu $ so that $A'=\psi^*\circ B$.
Then $B$ is $\nu$-linear, starts with the identity and
$$
B(u*''v)=Bu*'''Bv
$$
so that $B$ is an equivalence -- in the usual sense -- between $*''$ and
$*'''$. Hence \cite{GR}

\begin{prop}
\label{der:isos}
Any isomorphism between two differential star products 
on symplectic manifolds is the
combination of a change of parameter and a $\nu$-linear isomorphism. Any
$\nu$-linear isomorphism between two star products $*$ on $(M,\omega)$
and $*'$ on $(M',\omega')$ is the combination of the action on functions
of a symplectomorphism $\psi:M' \to M$ and an equivalence between $*$
and the pullback via $\psi$ of $*'$. In particular, it exists if and
only if those two star products are equivalent, i.e.\ if and only if
$(\psi^{-1})^*c(*')=c(*)$, where here $(\psi^{-1})^*$ denotes the action
on the second de~Rham cohomology space.
\end{prop}

In particular,
two differential star products $*$ on $(M,\omega)$ and 
$*'$ on $(M',\omega')$ are
isomorphic if and only if there exist $f(\nu)= \sum_{r\ge1} \nu^r f_r
\in \R\bbnu$ with $f_1\ne 0$ and $\psi \colon M' \to M$, a
symplectomorphism, such that $(\psi^{-1})^*c(*')(f(\nu))=c(*)(\nu)$. In
particular \cite{Guequiv}: if $\coh2=\R[\omega]$ then there is only
one star product up to equivalence and change of parameter.

 Omori et al.\ \cite{Maeda2} also show that
when reparametrizations are allowed then there is only one star product on
$\C P^n$.

A special case of Proposition \ref{der:isos} gives:

\begin{prop}
A symplectomorphism $\psi$ of a symplectic manifold  can be
extended to a $\nu$-linear automorphism of a given differential star
product on $(M,\omega)$ if and only if $(\psi)^*c(*)=c(*)$.
\end{prop}

Notice that this is always the case if $\psi$ can be connected to the
identity by a path of symplectomorphisms (and this result was in
Fedosov \cite{Fed2}).

\medskip
$\bullet$
If $G$ is a connected Lie group acting on the symplectic manifold $(M,\omega)$
by symplectomorphisms, each element of $G$ can be lifted to an automorphism
of a star product on $M$. The group $G$ acts on the quantum level if
there is a homomorphism $\rho$ from $G$ into the automorphism group of $*$
such that $\rho(g)=g^*+\nu\dots ~~\forall g\in G$.
At the Lie algebra level, one considers 
a homomorphism $\sigma$ from the  Lie algebra
$\g$ of $G$ to the algebra of derivations of $*$ such that 
$\sigma(X)=X^*+\nu\dots ~~\forall X\in \g$ where $X^*$ 
is the fundamental vector field on
$M$ associated to the action of $G$ on $M$ 
(i.e. $X^*_x=\frac{d}{dt} {\exp -tX }\cdot x_{\vert_0}).$
Now any local $\nu$-linear derivation of $*$ can be locally written
on a contractible set $U$ as 
$\ad_*(\frac{1}{\nu}f)$ for some $f\in C^\infty(U)[[\nu]]$.
\begin{definition}
Given a star product $*$ on a  Poisson manifold $(M,P)$
and given  a connected Lie group $G$ acting on $M$,
the star product gives {\bf { a quantization of the action 
of the Lie algebra $\g$}}
if there exists a {\bf{generalised moment map}} i.e. a map 
 $$
 \tau : \g \rightarrow \frac{1}{\nu}\CiM[[\nu]] : X \mapsto \tau(X)
 $$
 such that
 \begin{equation}\label{gemomap1}
 \tau(X)*\tau(Y)-\tau(Y)*\tau(X)=\tau([X,Y])\quad \forall X,Y\in\g
 \end{equation} 
and
 \begin{equation}\label{gemomap2}
 \tau(X)*u-u*\tau(X)=X^*u+\nu\dots ~~\forall X\in \g.
 \end{equation}
 \end{definition} 
Remark that the two conditions imply for the term of $\tau$ of
 order $-1$ in $\nu$ that the action of $G$ on $M$
admits a moment map, i.e.
there is a map  $\lambda: \g\rightarrow\CiM$ 
such that $\{\lambda(Y),u\}=Y^*u$ and 
$\{\lambda(Y),\lambda(Y)\}=\lambda(\{X,Y\})$,
and 
\[
\tau(X)=\frac{1}{\nu}\lambda(X) +\dots.
\]

When one can choose $\tau(X)=\frac{1}{\nu} \lambda(X)~\forall X\in\g$ the notion 
of a star product which is a quantization of the action 
of $\g$ coincides with the notion of a covariant star product.

When one can choose $\tau$ so that $\tau(X)*u-u*\tau(X)=X^*u~~\forall X\in \g$,
it implies that $X^*(u*v)=(X^*u)*v+u*(X^*v)$, so that 
\[
g^*(u*v)=g^*u*g^*v~\forall g\in G 
\]
which are the conditions in Bayen et al  \cite{BayenB}
for the star product to be {\bf geometrically invariant}; 
in that situation, the notion of  generalised moment map coincides
with the notion of quantum moment map introduced by  Xu in \cite{Xu}.

In the general case, when one has found a map $\tau$ at order $k$ in $\nu$
satisfying  condition (\ref{gemomap1})  at order $k$,
then one can extend things one order further if a  
Chevalley $2$-cocycle from $\g$
with values in $\CiM$ (for the representation of $Y\in\g$ given by $Y^*$)
is a $2$-coboundary. When one deals with geometrically 
invariant star products,
one can always find a $\tau$ such that 
$\tau(X)*u-u*\tau(X)=X^*u~~\forall X\in \g$
if $H^2(\g,\R)=0$ and $H^1(M;\R)=0$ as was found in \cite{Xu}. 

\section{Convergence of star products}

\begin{remark}
Let  $(M,\PO)$ be a Poisson manifold and let $*$ be a differential star
product on it with $1$ acting as the identity.
Observe that if there exists a value $k$ of $\nu$ such that
\[
u*v=\sum_{r=0}^\infty \nu^r C_r(u,v)
\]
converges (for the pointwize convergence of functions), for all $u,v\in\CiM$,
to $F_k(u,v)$ in such a way that $F_k$ is associative, then $F_k(u,v)=uv$.
\end{remark}
So assuming ``too much'' convergence kills all deformations. 
On the other hand,
in any physical situation, one needs some convergence properties to be able
 to compute the spectrum of quantum observables in terms of a star product
 as was done  for some observables already in Bayen, Flato, Fronsdal, 
 Lichnerowicz and Sternheimer \cite{BayenB2}.

\medskip
Now consider the example of Moyal star product on
 the symplectic vector space $(\R^{2n},\omega)$.
The formal formula for  Moyal star product
 \[
  (u*_Mv)(z)=
 \left.\exp \left(\frac{\nu}{2} \PO^{rs}
 \partial_{x^r}\partial_{y^s}\right)(u(x)v(y))\right\vert_{x=y=z} 
\]
obviously converges  when $u$ and $v$ are polynomials.

On the other hand, there is an integral formula for
 Moyal star product given by
\[
 (u*v)(\xi)=(\pi\hbar)^{-2n}\int u(\xi')v(\xi'') e^{\frac{2i}{\hbar}(\omega(\xi,\xi'')
 +\omega(\xi'',\xi')+\omega(\xi',\xi))}d\xi' d\xi'',
\]
and this product $*$ gives a structure of associative algebra 
on the space of rapidly
decreasing functions ${\CALS}(\R^{2n})$.

The formal formula  converges (for $\nu=i\hbar$) in the topology of ${\CALS}'$
 for $u$ and $v$ with compact Fourier transform.
 
 \medskip  With Michel Cahen and John Rawnsley, we  used 
 the method of quantization
of K\"ahler manifolds due to Berezin, \cite{Berezin}, as the inverse
of taking symbols of operators,   to
construct  on Hermitian symmetric spaces star products   
which are convergent on a large class of functions
on the manifold. I shall develop this construction in this 
last part of my talk.

\medskip Let me mention, before closing this introduction about convergence, 
the work of Rieffel \cite{Rieffel} where he introduces the notion of strict 
deformation quantization. An example of strict Fr\'echet quantization
has been recently given by Omori, Maeda, Niyazaki and Yoshioka in \cite{OMY3}.

Also very important  are the constructions of  operator representations 
of star products, in particular the works of Fedosov \cite{Fed2} and of
Bordemann, Neumaier and Waldmann
\cite{Bord1,Bord2}.

\subsection{Convergence of Berezin type star products on Hermitian 
symmetric spaces}
The method to construct a star product 
involves making a correspondence between operators and
functions (their Berezin symbols), transferring the operator
composition to the symbols, introducing a suitable parameter into
the Berezin composition of symbols, taking the asymptotic expansion in
this parameter on a large algebra of functions and then showing that
the coefficients of this expansion satisfy the cocycle conditions to
define a star product on the smooth functions.
The idea of an asymptotic expansion appeared already in Berezin 
\cite{Berezin} and  in  Moreno and Ortega-Navarro\cite{Moreno, Moreno2}.   

In  \cite{CGRII} we show that this asymptotic expansion
exists for compact $M$, and defines an associative multiplication on
formal power series in $k^{-1}$ with coefficients in $C^\infty(M)$
for compact coadjoint orbits. We also show that this formal power series
converges on the space of symbols for $M$ a Hermitian symmetric
space of compact type.

In \cite{Kara2}, Karabegov proves 
convergence for general compact coadjoint orbits (i.e. flag manifolds). 

In \cite{CGRIV} we study general Hermitian symmetric
spaces of non-compact type, and use their realisation as bounded
domains to define an analogous algebra of symbols of polynomial
differential operators.

Recently Reshetikhin and Takhtajan have announced \cite{RT} an associative
formal star product given by an asymptotic expansion on any K\"ahler manifold.
This they do in two steps, first building an associative product
for which $1$ is not a unit element, then passing to a star product. 

\subsubsection{Berezin symbols}
We denote by $(L,\nabla,h)$ a quantization bundle for the K\"ahler
manifold $(M,\omega,J)$ (i.e.
a holomorphic line bundle $L$ with connection $\nabla$ admitting an
invariant hermitian structure $h$, such that the curvature is curv$(\nabla)
= -2i\pi\omega$). We denote by $\H$ the Hilbert space of
square-integrable holomorphic sections of $L$ which we assume to be
non-trivial. The coherent states are vectors $e_q \in \H$ such that
$$
s(x) = \langle s, e_q \rangle q, \qquad \forall q \in \L_x,\quad x \in
M,\quad s \in \H
$$
where $\L$ denotes the complement of the zero-section in $L$. The
function
$$
\epsilon(x) = \vert q \vert^2 \Vert e_q \Vert^2, \qquad q \in
\L_x
$$
is well-defined and real analytic.

We introduce also the $2$-point function
$$
\psi(x,y) = {{\vert \langle e_{q^\prime}, e_q \rangle \vert^2} \over
{\Vert e_{q^\prime} \Vert^2 \Vert e_q \Vert^2}},\qquad q\in \L_x,\quad
q^\prime \in \L_y
$$
which is a globally defined real analytic function on $M\times M$
provided $\epsilon$ has no zeros. It is
a consequence of the Cauchy--Schwartz inequality that $\psi(x,y) \le 1$
everywhere, with equality where the lines spanned by $e_q$ and
$e_{q^\prime}$ coincide ($q\in \L_x$, $q^\prime \in \L_y$). 

Let $A:\H\to\H$ be a bounded linear operator and let
$$
\widehat{A}(x)={\langle Ae_q,e_q\rangle\over\langle e_q,e_q\rangle},\qquad
q\in \L_x,\quad x\in M
$$
be its symbol. The function $\widehat{A}$ has an analytic continuation to an
open neighbourhood of the diagonal in $M\times\overline M$ given by
$$
\widehat{A}(x,y) = {\langle Ae_{q^\prime}, e_q\rangle\over\langle
e_{q^\prime}, e_q\rangle}, \qquad q\in \L_x,\quad q^\prime \in \L_y
$$
which is holomorphic in $x$ and antiholomorphic in $y$.  We denote
by $\ehat{}$ the space of symbols of bounded operators on $\H$.
We can extend this definition of symbols to some unbounded operators
provided everything is well defined. 

\subsubsection{Composition of operators - Parameter}
The composition of operators on $\H$ gives rise to a product for the
corresponding symbols, which is associative and which we shall denote by
$*$ following Berezin, \cite{Berezin}. 
 The product  $*$ of symbols is given in terms of the
symbols by the integral formula
$$
(\AB )(x) = \int_M\widehat{A}(x,y)\widehat{B}(y,x)\psi(x,y)\epsilon(y)
{\omega^n(y)\over n!}.
$$
This formula is derived by use of the
adjoint $A^*$ of $A$ so 
 to apply it to the case where the operators
are unbounded we need to be able to use the adjoint of $A$ on coherent
states. To be able to take the symbol of the composition,
the result of applying $B$ to a coherent state
must be in the domain of $A$.

{\bf{Example: }}
 The identity map has symbol $\widehat{I} = 1$ and so $\widehat{A}*1 =
1*\widehat{A} = \widehat{A}$ for any operator $A$. In particular
$1*1=1$.

\medskip

Let $k$ be a positive integer. The bundle $(L^k =\otimes^kL,\nabla^k,
h^k)$ is a quantization bundle for $(M,k\omega ,J)$ and we denote by
${\H}^k$ the corresponding space of holomorphic sections and by
$\ehat{k}$ the space of symbols of linear operators on ${\H}^k$.
We let $\epsilonk$ be the corresponding function. We say that
{\bf the quantization is  regular} if $\epsilonk$ is a non-zero constant
for all nonnegative $k$ and if $\psi(x,y)=1$ implies $x=y$. The
significance of these conditions has been explained in \cite{CGRII}. 

 Let $G$  be a Lie group of isometries of the K\"ahler
manifold $(M,\omega,J)$ which lifts to a group of automorphisms of the
quantization bundle  $(L,\nabla,h)$. This automorphism group acts
naturally on  $(L^k,\nabla^{(k)},h^k)$; if  $g\in G$ and
 if $e^{(l)}_q$  is a coherent state of $L^l$, then 
$g.e^{(l)}_q=e^{(l)}_{gq}$ so the function $\epsilonk$
is invariant under $G$. In particular, if the quantization
is homogeneous, all $\epsilonk$ are constants.

In the regular case, 
the function $\psi$  in the integral  defining the composition of symbols
 for powers $L^k$ gets replaced
by powers $\psi^k$.

\medskip

 When  the manifold is compact, we have proven the following facts:

- when $\epsilon^{(k)}$ is constant for all $ k$  
 one has the nesting property $\ehat{k}\subset\hat
E(L^{k+1})$;

- with the same assumption $\cup_k\ehat{k}$ is dense in ${\mathcal{C}}
 ^{\circ}(M)$.
   
From the nesting property, one sees
that if $\hat A,\hat B$ belong to $\ehat{l}$ and if $k \geq l$ one may define
\begin{equation}\label{eqstark}
(\AkB)(x)=\int_M\hat A(x,y)\hat B(y,x)\psi^k(x,y)
\epsilon^{(k)}{k^n\omega^n\over
n!}.
\end{equation}
More generally, if $\widehat{A}$ is a symbol of an operator then its analytic
continuation $\widehat{A}(x,y)$ may have singularities where
$\psi(x,y)=0$ but $\widehat{A}\psi$ is always globally defined on
$M\times M$. If $M$ is not compact $\widehat{A}\psi$ may not be bounded,
so we introduce the class $\mathcal{B} \subset C^\infty(M)$ of functions $f$
which have an analytic continuation off the diagonal in
$M\times\overline{M}$ so that $f(x,y)\psi(x,y)^l$ is globally defined,
smooth and bounded on $K\times M$ and on $M\times K$ for each compact
subset $K$ of $M$ for some positive power $l$ and denote by ${\mathcal{B}}_l$
those for which the power $l$ suffices. Since $\psi$ is smooth and
bounded it is clear that ${\mathcal{B}}$ is a subalgebra of $C^\infty(M)$. In
the case $M$ is compact we obviously have $\ehat l \subset{\mathcal{B}}_l$.
If $\hat A,\hat B$ belong to $\mathcal{B}$, formula (\ref{eqstark}) is well
 defined for $k$ large enough.
 
  \medskip
 We study the behaviour of the integral \ref{eqstark} in terms of $k$.

\subsubsection{ An asymptotic formula}

In order to localise the integral \ref{eqstark} we use a version
of the Morse Lemma adapted from Combet, as in  Moreno and Ortega-Navarro. 

Let $(M,\omega,J)$ be a K\"ahler manifold with metric $g$. We denote by
$\exp_xX$ the exponential at $x$ of $X \in T_xM$. If $g$ is not complete
the exponential map may not be defined for all $x$ and $X$, but in any
case there is an open subset $V \subset TM$ where it is defined and
which contains the zero-section. The differential of the exponential map
at $0$ is the identity so the map $\alpha: V\to M \times M$ given by
$\alpha(X) = (p(X),\exp_{p(X)}X)$ where $p$ is the projection in the
tangent bundle $p \colon TM\to M$ is a diffeomorphism near the
zero-section. At any point of the zero-section the differential of
$\alpha$ is the identity.

\begin{prop}
Let $(M,\omega,J)$ be a K\"ahler manifold with metric $g$ and $\alpha
\colon V \to M\times M$ be the map defined above. Let $(L,\nabla,h)$ be
a regular quantization bundle over $M$  and let $\psi$ be the
corresponding $2$-point function on $M \times M$. Then there exists an
open neighbourhood  $W\subset V$ of the zero-section in $TM$ and a
smooth open embedding $\nu: W \to TM $ such that
$$
(-\log\psi{\scriptstyle\circ}\alpha{\scriptstyle\circ}\nu)(X)=
\pi g_{p(X)}(X,X),\quad X\in W
$$
and the differential of $\nu$ at any point of the zero-section is the
identity.
\end{prop}

Denote by $\widetilde{\mathcal{B}}$ the set of functions $f$ on $M\times M
\setminus \psi^{-1}(0)$ such that $f(x,y)\psi(x,y)^l$ has a smooth
extension to all of $M\times M$ which is bounded on $K\times M$ for each
compact subset $K\subset M$ for some $l$ and denote by 
${\widetilde{\mathcal{B}}}_l$ those for which the power $l$ suffices. 
If $f$, $g$ are in ${\mathcal{B}}$ then $f(x,y)g(y,x)$ is in 
$\widetilde{\mathcal{B}}$. Note also that if $f
\in \widetilde{\mathcal{B}}$ then its restriction to the diagonal
$\widehat{f}(x) = f(x,x)$ is smooth.

\medskip

For any $f$ belonging to ${\widetilde{\B}}_l$,
 the integral
$$
F_k(x)=\int_M f(x,y) \psi(x,y)^k k^n {{\omega^n(y)}\over{n!}}, \quad
\rm{for}\ k\geq l+1
$$
admits an asymptotic expansion
$$
F_k(x)\sim\sum_{r\geq 0}k^{-r}C_r(\hatf)(x)
$$
where $C_r$ is a smooth differential operator of order $2r$ depending
only on the geometry of $M$. The leading term is given by
$C_0(\hatf)(x) = \widehat{f}$.

We are not claiming that for this very general class of functions
$f(x,y)$ the integral  depends smoothly on $x$, only
that the coefficients of the asymptotic expansion do.
\medskip

In the regular case $\epsilonk$ has an asymptotic expansion $\sum_{r \ge
0} \epsilon_r /k^r$ as $k$ tends to infinity with $\epsilon_0=1$.
Indeed,
$$
1 = 1*_k1 = \epsilonk \int_M \psi(x,y)^k k^n {{\omega^n(y)}\over{n!}};
$$
 has an asymptotic expansion in $k^{-1}$ with leading term
$1$ by the previous proposition and we  can then invert the asymptotic
expansion to obtain one for $\epsilonk$.

\begin{thm}
Let $(M,\omega,J)$ be a K\"ahler manifold and $(L,\nabla,h)$  be a
regular quantization bundle over $M$. Let $\widehat{A},\widehat{B}$ be
in ${\mathcal{B}}$. Then
$$
(\AkB)(x)=\int_M\widehat{A}(x,y)\widehat{B}(y,x)\psi^k(x,y)\epsilonk k^n
{\omega^n(y)\over n!}(y),
$$
defined for $k$ sufficiently large, admits an asymptotic expansion in
$k^{-1}$ as $k\to\infty$
$$
(\AkB)(x)\sim \sum_{r\geq 0}k^{-r}C_r(\widehat{A},\widehat{B})(x)
$$
and the cochains $C_r$ are smooth bidifferential operators, invariant
under the automorphisms of the quantization and determined by the
geometry alone. Furthermore
$$
C_0(\widehat{A},\widehat{B})=\widehat{A}\widehat{B},
$$
and
$$
C_1(\widehat{A},\widehat{B}) -
C_1(\widehat{B},\widehat{A}) = {i\over\pi} \{
\widehat{A},\widehat{B}\}.
$$
\end{thm}

\subsubsection{A convergent star product for flag manifolds}

 We first would like to show that
the asymptotic expansion obtained above defines an associative formal star
product.  For this we assume that $(M,\omega,J)$ is a
flag manifold.
 Reshetikhin and Takhtajan have announced an analogous result for
general K\"ahler manifolds.

 Observe that if $G$  is a Lie group of isometries of the K\"ahler
manifold $(M,\omega,J)$ which lifts to a group of automorphisms of the
quantization bundle  $(L,\nabla,h)$ it acts
naturally on $(L^k,\nabla^{(k)},h^k)$ :
\[
g^*(\AkB)(x)=(g^*\hat A{ }*{ }g^*\hat B)(x)
\]
for any  $\hat A,\hat B$
in $\hat E(L^l)$ and any $k \geq l$.
Observe also that the bidifferential operators $C_r$ depend on
 the geometry alone thus are invariant under $G$.

\begin{lemma} Let $(M,\omega,J)$ be a flag manifold with $M =
G/K$ where $G$ is a compact simply-connected Lie group and $K$ the
centralizer of a torus.  Assume the geometric quantization conditions
are satisfied and let $(L,\nabla,h)$ be a quantization bundle over
$M$.  Let $\CL = \cup_k \ehat{k}$  be the union of the symbol spaces.
Then $\CL$ coincides with the space $E$ of vectors in $C^{\infty}(M)$
whose $G$-orbit is contained in a finite dimensional subspace.
\end{lemma}

In the case of a flag manifold as above, the group $G$ lifts to a group of
automorphisms of the quantization bundle $(L^k,\nabla^{(k)},h^k)$ 
hence the map $\ehat{l}\otimes\ehat{l}\to \Ci{M}$ given by $\hat A \otimes
\hat B\mapsto \AkB$ intertwines the action of $G$ and the bidifferential
operators $C_r$ are invariant under $G$.  
Thus,
if $\hat A$, $\hat B$ belong to $\ehat{l}$,
there exists an integer $a(l)$  such that $\AkB$ belongs to
$\ehat{a(l)}$  for all $k\ge l$, and such that $C_r(\hat A,\hat B)$
 belongs to $\ehat{a(l)}$ for every integer $r$.

Consider now the asymptotic development:
$$
\AkB = \sum_{r=0}^N k^{-r} C_r(\hat A,\hat B) + R_N(\hat A,\hat B,k)
$$
where
$$
\lim_{k\rightarrow\infty} k^N R_N(\hat A,\hat B,k) = 0.
$$
 The above tells us that $R_N(\hat A,\hat B,k)$ belongs
to $\ehat{a(l)}$ where $a(l)$ is independent of $k$. So, we can write
\begin{eqnarray}
&~&(\AkB)*_k\hat C = \sum_{r=0}^N k^{-r}C_r(\hat A,\hat B)*_k\hat C +
R_N(\hat A,\hat B)*_k\hat C\cr
&~&\quad\quad = \sum_{r,s=0}^N k^{-r-s}C_s(C_r(\hat A,\hat B),\hat C) \cr
&~&+\sum_{r=0}^N k^{-r} R_N(C_r(\hat A, \hat B),\hat C,k)+ 
R_N(\hat A,\hat B,k)*_k\hat C. 
\end{eqnarray}
The last two terms multiplied by $k^N$  tend to zero when $k$ 
tends to infinity.

\begin {thm}  The asymptotic expansion $\sum_{r\geq
0}k^{-r}C_r(u,v)$ yields a formal associative deformation of the usual
product of functions in $\CL$.
 It is  a formal star product which extends to all of $\Ci{M}$, using
uniform convergence.
\end{thm}

 We prove that $\AkB$ is a rational function of $k$ with no pole at
infinity, when the flag manifold is a hermitian symmetric space, by using
structure theory of these spaces.  

\begin{thm}
 Let $M$ be a compact hermitian symmetric space and let
$(L,\nabla,h)$ be a quantization bundle over M. Let $L^k = \otimes^k L$
and let ${\H}^k$ be the space of holomorphic sections of $L^k$. Let $\ehat{k}$
be the space of symbols of operators on ${\H}^k$. If $\hat A$, $\hat B$ belong
to $\ehat{l}$ and $k \geq l$, the product $\hat A *_k \hat B$ depends
rationally on k and has no pole at infinity,
hence the asymptotic expansion of $\AkB$ is convergent.
\end{thm}

This result is generalised by Karabegov \cite{Kara2}:

\begin{thm}
 For any generalised flag manifold the $*_k$ product
of two symbols is a rational function of $k$ without pole at infinity.
\end{thm}

\subsubsection{ Star product on bounded symmetric domains}
~

\medskip
\noindent{
\bf{Bounded  symmetric domains}} 

Let $\D$ denote a bounded symmetric domain. We shall use the
Harish-Chandra embedding to realise $\D$ as a bounded subset of its
Lie algebra of automorphisms. More precisely, if $\G$ is the connected
component of the group of holomorphic isometries then $\D$ is the
homogeneous space $\G/\K$ where $\G$ is a non-compact semi-simple Lie
group and $\K$ is a maximal compact subgroup. Let $\g$ be the Lie
algebra of $\G$, $\k$ the subalgebra corresponding with $\K$, $\gc$
and $\kc$ the complexifications and $\GC$, $\KC$ the corresponding
complex Lie groups containing $\G$ and $\K$. The complex structure on
$\D$ is determined by $\K$-invariant abelian subalgebras $\m_+$ and
$\m_-$ with
$$
\gc= \m_+ + \kc + \m_-, \qquad [\m_+,\m_-] \subset \kc,
\qquad\overline{\m_+} = \m_-
$$
where $\overline{\hbox{\vrule height1ex depth0pt width0pt}~~}$
denotes conjugation over the real form $\g$ of $\gc$.  The
exponential map sends $\m_{\pm}$ diffeomorphically onto subgroups
$M_{\pm}$ of $\GC$ such that $M_+\KC M_-$ is an open set in $\GC$
containing $\G$ and the multiplication map
$$
M_+\times\KC\times M_-\to M_+\KC M_-
$$
is a diffeomorphism. $\KC M_-$ is a parabolic
subgroup of $\GC$ and the quotient $\GC/\KC M_-$ a generalised flag
manifold. The $\G$-orbit of the identity coset can be identified
with $\G/\K$ and lies inside $M_+\KC M_-/\KC M_- \cong M_+$.
Composing this identification with the inverse of the exponential map
gives the desired Harish-Chandra embedding of $\D$ as a bounded open
subset of $\m_+$. We shall assume from now on that $\D \subset \m_+$ via
this embedding. In this realisation it is clear that the action of $\K$
on $\D$ coincides with the adjoint action of $\K$ on $\m_+$.

Following Satake, we define maps
$$
 k \colon \D\times\D \to \KC,\qquad
m_{\pm} \colon \D\times\D \to \m_{\pm}
$$
by
$$
\exp - \overline{Z'}\exp Z = \exp m_+(Z, Z') \,\, k(Z,
Z')^{-1} \exp m_-(Z, Z'), 
$$
for $Z, Z' \in \D.$ They satisfy
$$
\overline{ k(Z, Z')} =  k(Z', Z)^{-1}, \qquad 
\overline{m_+(Z, Z')} = -m_-(Z', Z).
$$

\medskip
\noindent{\bf{The holomorphic quantization}}

For any unitary character $\chi$ of $\K$ there is a Hermitian
holomorphic line bundle $L$ over $\D$ whose curvature is the
K\"ahler form of an invariant Hermitian metric on $\D$. If $\chi$
also denotes the holomorphic extension to $\KC$ then the Hermitian
metric has K\"ahler potential $\log\chi (k(Z,Z))$ and $L$ has a
zero-free holomorphic section $s_0$ with
$$
\vert s_0(Z) \vert^2 = \chi (k(Z,Z))^{-1}.
$$
For $\chi$ sufficiently positive $s_0$ is square-integrable and the
representation $U$ of $\G$ on the space $\H$ of square-integrable
sections of $L$ is one of Harish-Chandra's holomorphic discrete series.
 $s_0$ is a highest weight vector for the extremal
$K$-type so is a smooth vector for the representation. We form the
coherent states $e_q$ and see that $e_{s_0(0)}$ transforms the same
way as $s_0$ and so they must be equal up to a multiple. This means
that the coherent states are also smooth vectors of the
representation. Further, since the quantization is homogeneous,
$\epsilon$ will be constant. Thus
$$
\langle e_{s_0(Z')}, e_{s_0(Z)} \rangle = \epsilon\, \chi (k(Z,Z')).
$$

\begin{lemma}
Up to a constant (determined by the normalization of Haar measure on
$\G$) $\epsilon$ is the formal degree $d_U$ of the discrete series
representation $U$, hence by Harish-Chandra's formula,it
is a polynomial function of the differential $d\chi$ of
the character $\chi$.
\end{lemma}

The two-point function $\psi$ is given by
$$
\psi(Z,Z') = \frac{\vert \chi(k(Z,Z')) \vert^2}
{\chi(k(Z,Z))\chi(k(Z',Z'))};
$$
it takes the value 1 only on the
diagonal. 
\medskip

\noindent{\bf{Polynomial differential operators and symbols}}

We let $\A$ denote the algebra of holomorphic differential operators on
functions on $\D$ with polynomial coefficients. We filter $\A$ by both
the orders of the differentiation and the degrees of the coefficients:
$\A_{p,q}$ denotes the subspace of operators of order at most $p$ with
coefficients of degree at most $q$. Obviously, the composition of
operators gives a map
$$
\A_{p,q}\times\A_{p',q'} \to \A_{p+p',q+q'}.
$$

The global trivialization by $s_0$ of the holomorphic line bundle $L$
corresponding with the character $\chi$ allows us to transport the
above operators to act on sections of $L$ by sending $D \in \A$ to
$D^{\chi}$ where
$$
D^{\chi}(fs_0) = (Df)s_0.
$$
Let $\A(\chi)$ denote the resulting algebra of operators on sections of
$L$ and $\A_{p,q}(\chi)$ the corresponding subspaces.

In this non-compact situation elements of $\A(\chi)$ do not define
bounded operators on the Hilbert space $\H$, but the fact that the
coherent states are smooth vectors of the holomorphic discrete series
representation and that polynomials are bounded on $\D$ means that
each operator in $\A(\chi)$ maps the coherent states into $\H$ so that
it makes sense to speak of the symbols of these operators.

\begin{lemma}
The analytically continued symbol $\widehat {D^{\chi}} (Z,Z')$ of an
operator $D^{\chi}$ in $\A_{p,q}(\chi)$ is a polynomial in $Z$ and
$m_-(Z,Z')$ of bidegree $p,q$.

The space of symbols of the operators in $\A_{k,l}(\chi)$
 is the space
of polynomials in $Z$ and $m_-(Z,Z')$ of bidegree $k,l$ so, in
particular, is independent of $\chi$.
\end{lemma}

Denote by $\E_{p,q}$ the space of polynomials in $Z$ and $m_-(Z,Z')$ of
bidegree $p,q$ and by $\E$ the union of these spaces. $\E$ is an
algebra under pointwize multiplication, the algebra of symbols .
If we take a symbol in
$\E_{p,q}$ then it is the symbol of an operator in $\A_{p,q}(\chi)$.
Taking two such operators and composing them corresponds with the
composition of two polynomial operators in $\A_{p,q}$ and so can be
expressed in terms of a basis for $\A_{2p,2q}$ as a rational function
of $d\chi$. In other words the Berezin product $f * g$ of two symbols
$f$, $g$ in $\E_{p,q}$ is a symbol in $\E_{2p,2q}$ depending
rationally on $d\chi$.

\medskip

\noindent{\bf{The star product}}

We construct a
formal deformation of the algebra $C^\infty(\D)$ by first
constructing it on the subalgebra $\E$.

We consider the powers $L^k$ of the line bundle $L$ which correspond
with the powers $\chi^k$ of $\chi$. These powers have differentials
$kd\chi$, so the Berezin product $f*_kg$ of two symbols $f$, $g$ in
$\E_{p,q}$ is a rational function of $k$ by the results of the
previous section.

The symbols in $\E_{p,q}$ are in  $\B_0$
(it is enough to show that $\ad m_-(Z,Z')$
is bounded on $X\times\D$ and $\D\times X$ for any compact subset $X$
of $\D$ ); thus the
asymptotic expansion  exists.

 Since the Berezin product is associative for each $k$,
the same argument as in the compact situation shows that its asymptotic
expansion in $k^{-1}$ is an associative formal deformation on $\E$
with bidifferential operators as coefficients. To see that it extends
to all of $C^\infty(\D)$ we  show that $\E$ contains enough
functions to determine these operators, hence
the asymptotic expansion of $f*_kg$
has bidifferential operator coefficients which satisfy the cocycle
conditions to define a formal product on $C^\infty(\D)$ which is
associative.

\begin{thm}
Let $\D$ be a bounded symmetric domain and $\E$ the algebra of symbols
of polynomial differential operators on a homogeneous holomorphic line
bundle $L$ over $\D$ which gives a realisation of a holomorphic
discrete series representation of $\G$ (i.e 
$\E$ is the algebra of functions on $\D$
which are polynomials in $Z$ and $m_-(Z,Z)), 
$then for $f$ and $g$ in $\E$
the Berezin product $f*_kg$ has an asymptotic expansion in powers of
$k^{-1}$ which converges to a rational function of $k$. The
coefficients of the asymptotic expansion are bidifferential operators
which define an invariant and covariant star product on
$C^\infty(\D)$.
\end{thm}


\begin{thebibliography}{99}

\bibitem{Andler} M.~Andler, A.~Dvorsky and S.~Sahi, Kontsevich
 Quantization and
 invariant distributions on Lie groups, preprint math/9910104
 and math/9905065.

\bibitem{Arnal1} D.~Arnal, 
Le produit star de Kontsevich sur le dual d'une
alg\`ebre de Lie nilpotente.
\textit{C. R. Acad. Sci. Paris S\'er. I Math.}, 
{237} (1998) 823-826.


\bibitem{Arnal} D.~Arnal, N.~Ben Amar and M.~Masmoudi, Cohomology
of good graphs and Kontsevich linear star products, 
\textit{Lett. in Math. Phys.}
48 (1999) 291--306.

\bibitem{ACG} D.~Arnal, M.~Cahen and S.~Gutt ,
 Deformations on coadjoint orbits,
\textit {J. Geom. Phys.} 3 (1986) 327--351.

\bibitem{ArnCor}
D. Arnal, {$*$ products and representations of nilpotent Lie groups},
\textit{Pacific J. Math.} 114 (1984) 285--308 and
D.~Arnal and J.-C.~Cortet,
{$*$ products in the method of orbits for nilpotent Lie groups},
\textit{J. Geom. Phys.} 2 (1985) {83--116}

\bibitem{ArnCor2} D.~Arnal and J.-C.~Cortet, 
{Nilpotent Fourier-transform and applications},
\textit{Lett. Math. Phys.} {9} (1985) 25--34 and  D. Arnal and S. Gutt,
{D\'ecomposition de $L\sp 2(G)$ et transformation de 
Fourier adapt\'ee pour un groupe $G$ nilpotent},
\textit{C. R. Acad. Sci. Paris S\'er. I Math.} 306 (1988) 25--28.

\bibitem{ACMP} D.~Arnal, J.-C.~Cortet, P.~Molin and G.~Pinczon, 
Covariance and geometrical invariance in star quantization, 
\textit{Journ. of Math. Phys.} 24 (1983) 276--283.

\bibitem{AL} D.~Arnal, J.~Ludwig and M. Masmoudi, 
D\'eformations covariantes sur les
orbites polaris\'ees d'un groupe de Lie, 
\textit{ Journ. of Geom. and Phys.} 14 (1994) 309--331.
 
\bibitem{Asin} S.~Asin, PhD thesis, Warwick University 1998.

\bibitem{BayenB} F.~Bayen, M.~Flato, C.~Fronsdal,
A.~Lichnerowicz and D.~Sternheimer,
{Quantum mechanics as a deformation of classical mechanics},
\textit{Lett. Math. Phys.} 
1 (1977) 521--530 and 
Deformation theory and quantization, part I,
\textit{Ann. of Phys.} 111 (1978) 61--110. 

\bibitem{BayenB2} F.~Bayen, M.~Flato, C.~Fronsdal,
A.~Lichnerowicz and D.~Sternheimer,
Deformation theory and quantization, part II,
\textit{Ann. of Phys.} 111 (1978) 111--151  


\bibitem{Berezin} F.A. Berezin, {General concept of quantization},
 \textit{Commun. Math. Phys.} {40} {(1975)} {153--174}.

\bibitem{Bertelson} M.~Bertelson,
Equivalence de produits star,
\textit{M\'emoire de Licence} U.L.B. (1995)
and  M.~Bertelson, M.~Cahen and S.~Gutt,
Equivalence of star products,
\textit{Class. Quan. Grav.} 14 (1997) A93--A107.

\bibitem{BBG} M.~Bertelson, P.~Bieliavsky and S.~Gutt, Parametrizing
equivalence  classes of invariant star products, 
\textit{Lett. in Math. Phys.} 46
 (1998) 339--345.
 
\bibitem{BP} F. Bidegain, G. Pinczon, Quantization of Poisson-Lie groups
 and applications, \textit{Commun. Math. Phys.} 179 (1996) 295--332.
 
 \bibitem{BP2} F. Bidegain, G. Pinczon, A $*$-product approach
 to non-compact quantum groups, \textit{Lett. Math. Phys.} 33 (1995) 231--240.

 
\bibitem{BFGP} P. Bonneau, M. Flato, M. Gerstenhaber, G. Pinczon,
The hidden group structure of quantum groups: strong duality, rigidity and
preferred deformations, \textit{Commun. Math. Phys.} 161 (1994) 125--156.

\bibitem{Bord1} M. Bordemann, N. Neumaier and S. Waldmann, Homogeneous
Fedosov star products on cotangent bundles I, \textit{ Comm. in Math. Phys.}
198 (1998) 363--396.

\bibitem{Bord2} M. Bordemann, N. Neumaier and S. Waldmann, Homogeneous
Fedosov star products on cotangent bundles II, 
\textit{ Journ. of Geom. and Phys.}
29 (1999) 199--234.

\bibitem{CDG} M.~Cahen, M. De Wilde and S.~Gutt, Local cohomology of the algebra
of smooth functions on a connected manifold, \textit{Lett. in Math. Phys.} 4 
 (1980) 157--167.


\bibitem{CFGS} M.~Cahen, M.~Flato, S.~Gutt and D.~Sternheimer,
Do different deformations lead to the same spectrum~?, 
\textit{Journ. of Geom. and 
		    Phys.} 2 (1985) 35--48.

\bibitem{CGparalman} M.~Cahen and S.~Gutt, Regular $*$ representations 
of Lie Algebras, \textit{Lett. in Math. Phys.}  6 (1982) 395--404.

\bibitem{CG} M.~Cahen and S.~Gutt, Produits $*$ sur les 
		orbites des groupes semi-simples 
		    de rang 1, \textit{C.R. Acad. Sc. Paris}  296 (1983) 821--823 and
          An algebraic construction of $*$ product on the regular
		    orbits of semisimple Lie groups, \textit{Bibliopolis Ed. 
		    Naples, Volume in honour 
		    of I. Robinson} (1987) 71--82 .

		    
\bibitem{CaGulocsym} M.~Cahen and S.~Gutt, Produits $*$ sur les espaces 
affins symplectiques localement sym\'etriques'', 
\textit{C.R. Acad. Sc. Paris} 297 (1983)
 417--420.    

\bibitem {CGRII}  M.~Cahen, S.~Gutt and J.~Rawnsley, Quantisation of  
K\"ahler manifolds II, \textit{Transactions A.M.S.} 337 (1993) 73--98.	

\bibitem {CGRIV}  M.~Cahen, S.~Gutt and J.~Rawnsley,
Quantisation of K\"ahler manifolds III and 
IV, \textit{Lett. in Math. Phys.} 30 
 (1994) 291--305 and {Lett. in  Math. Phys.} 34 
 (1995) 159--168.
 
\bibitem{CGR} M.~Cahen, S.~Gutt and J.~Rawnsley, 
  On tangential star products for the coadjoint Poisson structure, 
 \textit{Comm. in Math. Phys.} 180 (1996) 99--108.
 
\bibitem{Chloup} V.~Chloup, Star products on the algebra of
polynomials on the dual of a semi-simple Lie
algebra, \textit{Acad. Roy. Belg. Bull. Cl. Sci.}  8 (1997)
263--269.

\bibitem{Connes} A.~Connes, Non commutative differential geometry,
IHES Publ. Math. 62 (1985) 257--360.

\bibitem{CFS} {A. Connes, M. Flato and D. Sternheimer},
{Closed star products and cyclic cohomology},
 \textit{Lett. Math. Phys.} {24} (1992) {1--12}.


\bibitem{Deligne} P.~Deligne,
D\'eformations de l'Alg\`ebre des Fonctions d'une Vari\'et\'e
Symplectique: Comparaison entre Fedosov et De~Wilde Lecomte, 
\textit{Selecta Math.~(New series).} 1 (1995) 667--697.

\bibitem{DeWilde} M.~De~Wilde,
Deformations of the algebra of functions on a symplectic manifold:
a simple cohomological approach. Publication no. 96.005, Institut
de Math\'ematique, Universit\'e de Li\`ege, 1996.

\bibitem{DL} M.~De~Wilde and P.~Lecomte,
Existence of star-products and of formal deformations of the Poisson Lie
algebra of arbitrary symplectic manifolds,
\textit{Lett. Math. Phys.} 7 (1983) 487--496.

\bibitem{DL2} M.~De~Wilde and P.~Lecomte,
 Formal deformations of the Poisson Lie algebra of a symplectic
manifold and star products: existence, equivalence, derivations, 
\textit{ in Deformation Theory of Algebras and Structures and Applications},
ed. by Hazewinkel and Gerstenhaber, Kluwer (1988) 897--960.

\bibitem{DGL} M. De~Wilde, S. Gutt and P.B.A.~Lecomte,
\`A propos des deuxi\`eme et troisi\`eme espaces de cohomologie de
l'alg\`ebre de Lie de Poisson d'une vari\'et\'e symplectique.
\textit{Ann. Inst. H. Poincar\'e Sect. A (N.S.)} 40 (1984) 77--83.

\bibitem{Dito} G.  Dito, 
Kontsevich star product on the dual of a Lie algebra,
\textit{ Lett. in Math. Phys.} 48 (1999) 307--322.

\bibitem{Drinfeld} V.G. Drinfeld, {Quantum Groups},
\textit{Proc. ICM86, Berkeley, Amer. Math. Soc.}
  1 (1987) {101--110}.

\bibitem{EtKaz} P.~Etingof and D.~Kazhdan, 
{Quantization of Lie Bialgebras I},
  \textit{ Selecta Math.}, new series 2 (1996) 1--41.
  
\bibitem{EtKaz2} P.~Etingof and D.~Kazhdan, 
{Quantization of Poisson algebraic groups and Poisson homogeneous spaces}, 
in A.~Connes et al (eds.) \textit{Sym\'etries quantiques} 
(Les Houches, 1995), North-Holland,
Amsterdam, (1998) 935--946 (also q-alg/9510020).


\bibitem{Fed} B.V.~Fedosov,
A simple geometrical construction of deformation quantization,
\textit{J. Diff. Geom.} 40 (1994) 213--238.

\bibitem{Fed2} B.V.~Fedosov,
\textit{Deformation quantization and index theory}.
Mathematical Topics Vol.\ 9, Akademie Verlag, Berlin, 1996.

\bibitem{Fedosov} B.V.~Fedosov, 
The index theorem for deformation quantization, in M.~Demuth et al. (eds.)
Boundary value problems, Schr\"odinger operators, deformation
quantization, \textit{Mathematical Topics} Vol.\ 8, Akademie Verlag,
 Berlin, (1996) 206--318.

\bibitem{Fed3} B.V.~Fedosov, 
On $G$-Trace and $G$-Index in deformation quantization,
preprint 99/31, Universit\"at Potsdam.

\bibitem{FL} R.~Fioresi, M.~A.~Lledo, On the deformation quantization of
coadjoint orbits of semisimple groups, preprint math/9906104.

\bibitem{Flato} M.~Flato, Deformation view of physical theories,
\textit{Czec. J. Phys.} B32 (1982) 472--475.
 
\bibitem{FLS1} M.~Flato, A.~Lichnerowicz and D.~Sternheimer,
 {D\'eformations $1$-diff\'erentiables d'alg\`ebres de Lie attach\'ees
  \`a une vari\'et\'e symplectique ou de contact},
  \textit{C. R. Acad. Sci. Paris S\'er. A} 279 (1974) 877--881 and
  \textit{Compositio Math.} 31 (1975) 47--82.
  
\bibitem{FLS} M.~Flato, A.~Lichnerowicz and D.~Sternheimer,
 Crochet de Moyal--Vey et quantification,
 \textit{C. R. Acad. Sci. Paris I Math.} 283 (1976)  19--24.

\bibitem{Fronsdal} C.~Fronsdal, Some ideas about quantization,
\textit{Reports On Math. Phys.} 15 (1978) 111--145.

\bibitem{Gerst} M.~Gerstenhaber, On the deformation of rings and algebras.
\textit{Ann. Math.} 79 (1964) 59--103.

\bibitem{Guequiv} S.~Gutt, Equivalence of deformations and 
 associated $*$ products, \textit{Lett. in Math. Phys.} 3 (1979) 297--309.
		    
\bibitem{Gutt23} S.~Gutt,
 Second et troisi\`eme espaces de cohomologie diff\'erentiable de
l'alg\`ebre de Lie de Poisson d'une vari\'et\'e symplectique,
\textit{Ann. Inst. H. Poincar\'e Sect. A (N.S.)} 33 (1980) 1--31.

\bibitem{Gutg} S.~Gutt, An explicit $*$ product on the 
cotangent bundle of a Lie group, 
\textit{Lett. in Math . Phys.} 7 (1983), 249--258.

\bibitem{Gutt} S.~Gutt, 
On some second Hochschild cohomology spaces for algebras
of functions on a manifold, \textit{ Lett. Math. Phys.} 39 (1997) 157--162.

\bibitem{GR} S.~Gutt and J.~Rawnsley,
 Equivalence of star products on a symplectic manifold;
an introduction to Deligne's \cech\ cohomology classes, 
\textit{ Journ.  Geom.  Phys.}
29 (1999) 347--392.

\bibitem{Kara2} A.~Karabegov, 
Berezin's quantization on flag manifolds and spherical
modules, \textit{Trans. Amer. Math. Soc.}  359 (1998) 1467--1479.

\bibitem{Karaequivsepvar} A.~Karabegov,
Cohomological classification of deformation quantisations
with separation of variables,
\textit{Lett. Math. Phys.} 43 (1998) 347--357.

\bibitem{Kara3} A.~Karabegov,
 On the canonical normalisation of a trace
 density of deformation quantization, 
 \textit{Lett. in Math. Phys.}  45 (1999) 217--228.


\bibitem{Kathotia} V.~Kathotia, 
Kontsevich universal formula for deformation
 quantization and the CBH formula, 
 preprint math/9811174.

\bibitem{K} M.~Kontsevich,
Deformation quantization of Poisson manifolds, I. IHES
preprint q-alg/9709040.

\bibitem{Lecomte} P.B.A. Lecomte,
  Application of the cohomology of graded Lie algebras to formal
   deformations of Lie algebras,
 \textit{Lett. Math. Phys.} 13 (1987) 157--166.
 
 \bibitem{Lich} A.~Lichnerowicz, Cohomologie $1$-diff\'erentiable
 des alg\`ebres de Lie attach\'ees \`a une vari\'et\'e symplectique 
 ou de contact, \textit{Journ. Math. pures et appl.} 53 (1974) 459--484.
  
\bibitem{Lichneequiv} A.~Lichnerowicz,  
 Existence and equivalence of twisted products on a symplectic manifold,
\textit{Lett. Math. Phys.} {3} (1979) 495--502.


\bibitem{Lichne} A.~Lichnerowicz,
D\'eformations d'alg\`ebres associ\'ees \`a une vari\'et\'e
symplectique (les $*_\nu$-produits),
\textit{Ann. Inst. Fourier, Grenoble} 32 (1982) 157--209.

\bibitem{Masmoudi} M.~Masmoudi, 
Tangential formal deformations of the Poisson 
bracket and tangential star products on a regular Poisson manifold,
 \textit{J. Geom. Phys.} 9 (1992) 155--171.
  
  \bibitem{Moreno} C. Moreno and P. Ortega-Navarro, 
$*$-products on $D^1(C)$, $S^2$ and related spectral analysis, 
\textit{Lett. Math. Phys.}
 7 {(1983)} 181--193.
 
 \bibitem{Moreno2} C. Moreno,
Star-products on some K\"ahler-manifolds,
\textit{Lett. Math. Phys.}
 11 (1986) 361--372.

  \bibitem{Nadaud} F.~Nadaud, On continuous and differential Hochschild
cohomology, \textit{Lett. in Math. Phys.} 47 (1999) 85--95.

\bibitem{NeroVla} O.M.~Neroslavsky and A.T.~Vlassov,
Sur les d\'eformations de l'alg\`ebre des fonctions d'une vari\'et\'e
symplectique,
\textit{C. R. Acad. Sci. Paris S\'er. I Math.} 292 (1981) 71--76.

\bibitem{Nest-Tsygan} R.~Nest and B.~Tsygan,
Algebraic index theorem for families,
 \textit{Advances in Math.} 113 (1995) 151--205.

\bibitem{NT2} R.~Nest and B.~Tsygan, Algebraic index theorem,
 \textit{Comm. in Math. Phys.} 172 (1995) 223--262. 

\bibitem{Neumaier} N.~Neumaier, 
Local $\nu$-Euler Derivations and Deligne's Characteristic
Class of Fedosov Star Products and Star Products of Special Type, preprint 
math/9905176.

\bibitem{OMY1} H.~Omori, Y.~Maeda and A.~Yoshioka,
Weyl manifolds and deformation quantization,
 \textit{Adv. Math.} 85 (1991) 224--255.


\bibitem{Maeda2} H.~Omori, Y.~Maeda and A.~Yoshioka,
The uniqueness of star-products on ${\rm P}\sb n({\bf C})$,
in C.~H.~Gu et al. (eds.)
\textit{Differential geometry (Shanghai, 1991)}. pp 170--176.
World Sci. Publishing, River Edge, NJ, 1993.

\bibitem{OMYclosed} H. Omori and Y. Maeda and A. Yoshioka,
 Existence of a closed star product,
\textit{Lett. Math. Phys.} 26 (1992) {285--294}.


\bibitem{OMY2} H.~Omori, Y.~Maeda and A.~Yoshioka, Deformation
quantizations of Poisson algebras, in Y.~Maeda et al. (eds.),
symplectic geometry and quantization (Sanda and Yokohama, 1993) 
\textit{Contemp. Math.} 179 (1994) 213--240.

\bibitem{OMY3}  H.~Omori, Y.~Maeda, N. Niyazaki and A.~Yoshioka,
An example of strict Fr\' echet deformation quantization, preprint 1999.

\bibitem{Pinczon} G.~Pinczon, 
On the equivalence between continuous and differential
deformation theories, \textit{ Lett. Math. Phys.} 39 (1997) 143--156.

\bibitem{Rauch} D.~Rauch,
Equivalence de produits star et classes de Deligne,
\textit{M\'emoire de Licence} U.L.B. (1998).

\bibitem{CGRI} J.~Rawnsley, M.~Cahen and S.~Gutt, Quantization of 
		K\"ahler manifolds I,
		  \textit{Journal of Geometry and Physics} 7 (1990) 45--62.
		  
\bibitem{RT} N.~Reshetikhin and L.~Takhtajan, Deformation quantization of
K\"ahler manifolds, preprint math/9907171.

\bibitem{Rieffel} M. Rieffel, Questions on quantization, 
in L.~Ge et al. (eds.), Operator algebras and operator theory
(Shanghai,1997),\textit{ Contem. Math.}
228 (1998) 315--328.

\bibitem{Stern1} D. Sternheimer, Phase-space representations,
in M.~Flato et al. (eds.), Applications of group theory in physics and 
mathematical physics (Chicago, 1982), \textit{ Lect. in Appl. Math.} 21, 
Amer. Math. Soc., Providence RI,   (1985) 255-267.

\bibitem{Stern2} D. Sternheimer, 
Deformation Quantization  Twenty Years after,
in J.~RembieliÕnski (ed.), Particles, fields and gravitation (Lodz 1998)
\textit{ AIP conference proceedings} 453 (1998) 107--145.
and math/9809056.
 
\bibitem{Tam1} D. Tamarkin, Quantization of Poisson structures
on $\R^2$, preprint math/9705007.

\bibitem{Tam2} D. Tamarkin, Another proof of M. Kontsevich
formality theorem, preprint math/9803025, and
Formality of chain operad of small squares, preprint math/9809164.

\bibitem{Vey} J.~Vey,
D\'eformation du crochet de Poisson sur une vari\'et\'e
symplectique,
\textit{Comment. Math. Helvet.} 50
(1975) 421--454.

\bibitem{W} {A. Weinstein}, {Deformation quantization},
S\'eminaire Bourbaki 95,
\textit{Ast\'erisque} {227} (1995) 389--409.
 
\bibitem{WX} A.~Weinstein and P.~Xu,
Hochschild cohomology and characteristic classes for star-products, 
preprint q-alg/9709043.

\bibitem{Xu} Ping Xu, Fedosov $*$-products and quantum moment maps, 
\textit{Comm.
in Math. Phys.} 197 (1998) 167--197.


\end{thebibliography}
\end{document}